\documentclass[12pt,twoside,british]{article}
\usepackage[T1]{fontenc}
\usepackage[utf8]{inputenc}
\usepackage[a4paper]{geometry}
\usepackage{color}
\usepackage{mathrsfs}
\usepackage{amsmath}
\usepackage{amsthm}
\usepackage{amssymb}
\usepackage{setspace}

\makeatletter
\numberwithin{equation}{section}
\numberwithin{figure}{section}
\theoremstyle{plain}
\newtheorem{thm}{\protect\theoremname}[section]
\theoremstyle{plain}
\newtheorem{cor}[thm]{\protect\corollaryname}
\theoremstyle{plain}
\newtheorem{lem}[thm]{\protect\lemmaname}
\theoremstyle{definition}
\newtheorem{defn}[thm]{\protect\definitionname}
\theoremstyle{plain}
\newtheorem{prop}[thm]{\protect\propositionname}
\newcommand{\lyxaddress}[1]{
	\par {\raggedright #1
	\vspace{1.4em}
	\noindent\par}
}

\@ifundefined{date}{}{\date{}}
\usepackage[english]{babel}
\usepackage{latexsym}
\usepackage{amsthm}
\usepackage{eucal}
\usepackage{eufrak}
\usepackage{epstopdf}
\usepackage{graphicx}
\theoremstyle{plain}

\newtheoremstyle{boldremark}
    {\dimexpr\topsep/2\relax} 
    {\dimexpr\topsep/2\relax} 
    {}          
    {}          
    {\bfseries} 
    {.}         
    {.5em}      
    {}          

\theoremstyle{boldremark}
\newtheorem{brem} [thm] {Remark} 

\usepackage{fancyhdr}
\usepackage{blindtext}
\usepackage{titlesec}

\titleformat{name=\chapter}[display]
{\normalfont\Huge\itshape}
{\titlerule[1pt]\vspace{-33pt}\filleft%
  \parbox[t]{6em}{%
    \raggedleft%
    \rule{\linewidth}{0.5ex}\newline%
    \chaptertitlename\ \thechapter%
  }%
}
{4pc}
{\normalfont\upshape\bfseries\Huge}
\makeatother

\usepackage{babel}
\providecommand{\corollaryname}{Corollary}
\providecommand{\definitionname}{Definition}
\providecommand{\lemmaname}{Lemma}
\providecommand{\propositionname}{Proposition}
\providecommand{\theoremname}{Theorem}

\begin{document}
\begin{singlespace}

\title{\noindent \textbf{Besov regularity for a class of singular\\or degenerate
elliptic equations}}
\end{singlespace}
\begin{singlespace}

\author{\noindent Pasquale Ambrosio}
\end{singlespace}
\begin{singlespace}

\date{\noindent September 1, 2021}
\end{singlespace}
\maketitle
\begin{abstract}
\begin{singlespace}
\noindent Motivated by applications to congested traffic problems,
we establish higher integrability results for the gradient of local
weak solutions to the strongly degenerate or singular elliptic PDE
\[
-\mathrm{div}\left((\vert\nabla u\vert-1)_{+}^{q-1}\frac{\nabla u}{\vert\nabla u\vert}\right)=f,\,\,\,\,\,\mathrm{in}\,\,\Omega,
\]
where $\Omega$ is a bounded domain in $\mathbb{R}^{n}$ for $n\geq2$,
$1<q<\infty$ and $\left(\,\cdot\,\right)_{+}$ stands for the positive
part. We assume that the datum $f$ belongs to a suitable Sobolev
or Besov space. The main novelty here is that we deal with the case
of \textit{subquadratic growth}, i.e. $1<q<2$, which has so far been
neglected. In the latter case, we also prove the higher fractional
differentiability of the solution to a variational problem, which
is characterized by the above equation. For the sake of completeness,
we finally give a Besov regularity result also in the case $q\geq2$.\vspace{0.2cm}
\end{singlespace}
\end{abstract}
\begin{singlespace}
\end{singlespace}
\begin{singlespace}
\end{singlespace}
\begin{singlespace}

\section{Introduction}
\end{singlespace}

\begin{singlespace}
\noindent $\hspace*{1em}$Let us consider the following pair of variational
problems in duality: 
\begin{equation}
\inf_{\sigma\in L^{p}\left(\Omega,\mathbb{R}^{n}\right)}\left\{ \int_{\Omega}\mathcal{H}(\sigma(x))\,dx:-\mathrm{div}\,\sigma=f,\,\,\langle\sigma,\nu\rangle=0\,\,\mathrm{on}\,\,\partial\Omega\right\} \tag{\ensuremath{\mathrm{P1}}}\label{eq:P1}
\end{equation}
and
\begin{equation}
\sup_{u\in W^{1,p'}\left(\Omega\right)}\left\{ \int_{\Omega}u(x)f(x)\,dx-\int_{\Omega}\mathcal{H}^{*}(\nabla u(x))\,dx\right\} ,\tag{\ensuremath{\mathrm{P2}}}\label{eq:P2}
\end{equation}
where $\Omega$ is a bounded connected open subset of $\mathbb{R}^{n}$
($n\geq2$) with Lipschitz boundary, $f$ is a given function defined
over $\overline{\Omega}$ with zero mean (i.e. $\int_{\Omega}f\,dx=0$)
and the function $\mathcal{H}$ is defined for all $\sigma\in\mathbb{R}^{n}$
by
\begin{equation}
\mathcal{H}(\sigma):=\frac{1}{p}\left|\sigma\right|^{p}+\left|\sigma\right|,\label{eq:choice}
\end{equation}
where $p\in\left(1,+\infty\right)$. With such a choice, we get
\[
\nabla\mathcal{H}^{*}(z)=\left(\left|z\right|-1\right)_{+}^{q-1}\frac{z}{\left|z\right|},\,\,\,\,z\in\mathbb{R}^{n},
\]
where $\mathcal{H}^{*}$ is the Legendre transform of $\mathcal{H}$,
$q=p'=p/(p-1)$ is the conjugate exponent of $p$ and $\left(\,\cdot\,\right)_{+}$
stands for the positive part (see \cite{Br}). It is well known that
the Euler-Lagrange equation of the functional in (\ref{eq:P2}) is
given by the PDE
\begin{equation}
\begin{cases}
\begin{array}{cc}
-\mathrm{div}\,(\nabla\mathcal{H}^{*}(\nabla u))=f & \mathrm{in}\,\,\Omega,\\
\,\,\,\,\,\,\,\langle\nabla\mathcal{H}^{*}(\nabla u),\nu\rangle=0 & \mathrm{on}\,\,\partial\Omega,
\end{array}\end{cases}\label{eq:Eul-Lag}
\end{equation}
which has to be meant in the distributional sense.\\
$\hspace*{1em}$In this paper we prove local higher integrability
results both for the solution to (\ref{eq:P1}) and for the gradient
of the local weak solutions to the equation  appearing in (\ref{eq:Eul-Lag}),
when the datum $f$ belongs either to the Sobolev class $W_{loc}^{1,p}\left(\Omega\right)$
or to the local Besov space $B_{p,\infty,loc}^{\alpha}\left(\Omega\right)$,
for some $\alpha\in\left(0,1\right)$. Furthermore, we establish the
higher fractional differentiability of the solution to (\ref{eq:P1})
in the case $1<q<2$, under the same assumptions about the datum $f$.\\
$\hspace*{1em}$From the very beginning we warn the reader that, for
the sake of simplicity, we will confine ourselves to the case of the
cost function $\mathcal{H}$ defined in (\ref{eq:choice}), rather
than considering the more general function
\begin{equation}
\sigma\in\mathbb{R}^{n}\mapsto\frac{1}{p}\left|\sigma\right|^{p}+\delta\left|\sigma\right|,\label{eq:general}
\end{equation}
which leads to the equation
\[
-\mathrm{div}\left((\vert\nabla u\vert-\delta)_{+}^{q-1}\frac{\nabla u}{\vert\nabla u\vert}\right)=f,
\]
where $\delta>0$. Indeed, there is no loss of generality in supposing
$\delta=1$ and our proofs still work in the more general case (\ref{eq:general}),
with minor modifications.\\
$\hspace*{1em}$If we suppose that $f\in L^{p}\left(\Omega\right)$
and that the infimum in (\ref{eq:P1}) is finite, then problem (\ref{eq:P1})
consists in minimizing a strictly convex and coercive functional on
$L^{p}$ subject to a convex and closed constraint: therefore, a solution
$\sigma_{0}$ to (\ref{eq:P1}) exists and must be unique. Moreover,
we recall that by standard convex duality (see \cite{Eke} for example),
the values $\inf\,$(\ref{eq:P1}) and $\sup\,$(\ref{eq:P2}) coincide
and the primal-dual optimality condition characterizes the minimizer
$\sigma_{0}$ of (\ref{eq:P1}) through the equality 
\begin{equation}
\sigma_{0}(x)=\nabla\mathcal{H}^{*}(\nabla u_{0}(x)),\,\,\,\,\,\mathrm{for}\,\,\,\mathscr{L}^{n}\textrm{-}\mathrm{a.e.}\,\,\,x\in\Omega,\label{eq:opt-cond}
\end{equation}
where $u_{0}$ is a solution of (\ref{eq:P2}) and $\mathscr{L}^{n}$
denotes the $n$-dimensional Lebesgue measure. This is equivalent
to the requirement that $u_{0}$ is a weak solution of the Euler-Lagrange
equation (\ref{eq:Eul-Lag}), in the sense that
\[
\int_{\Omega}\langle\nabla\mathcal{H}^{*}(\nabla u_{0}(x)),\nabla\varphi(x)\rangle\,dx=\int_{\Omega}f(x)\varphi(x)\,dx,\,\,\,\,\,\mathrm{for}\,\,\mathrm{every}\,\,\varphi\in W^{1,q}\left(\Omega\right).
\]
Furthermore, since $f$ has zero mean, using the direct methods of
the Calculus of Variations it is not difficult to show that the dual
problem (\ref{eq:P2}) admits at least one solution $u_{0}$ belonging
to $W_{\diamond}^{1,q}\left(\Omega\right)$, where
\[
W_{\diamond}^{1,q}\left(\Omega\right):=\left\{ u\in W^{1,q}\left(\Omega\right):\,\int_{\Omega}u(x)\,dx=0\right\} ,
\]
and such that (\ref{eq:opt-cond}) holds, so that $u_{0}$ is a distributional
solution of the strongly degenerate or singular elliptic equation
\begin{equation}
-\mathrm{div}\left((\vert\nabla u\vert-1)_{+}^{q-1}\frac{\nabla u}{\vert\nabla u\vert}\right)=f,\label{eq:elliptic}
\end{equation}
under homogeneous Neumann boundary conditions. We also note that,
in general, if one looks at the solutions $u$ of the above equation,
no more than $C^{0,1}$ regularity should be expected for them: indeed,
every $1$--Lipschitz function $u$ is a solution of the homogeneous
equation. Moreover, when $q\geq2$ we have 
\[
\frac{(\left|\nabla u\right|-1)_{+}^{q-1}}{\left|\nabla u\right|}\left|\xi\right|^{2}\leq\langle D^{2}\mathcal{H}^{*}(\nabla u)\xi,\xi\rangle\leq(q-1)(\left|\nabla u\right|-1)_{+}^{q-2}\left|\xi\right|^{2},\,\,\,\,\xi\in\mathbb{R}^{n},
\]
that is, the Hessian matrix $D^{2}\mathcal{H}^{*}$ has eigenvalues
becoming degenerate in the region $\left\{ \vert\nabla u\vert\leq1\right\} $.
\\
$\hspace*{1em}$The main novelty of this paper is that we shall address
the singular case $1<q<2$, which has so far been neglected, since
extra technical difficulties arise concerning elliptic regularity
whenever we are in this case. This phenomenon, which also occurs in
the classical $q$-harmonic setting, has been very well explained
in \cite[Section 2.6]{Balci}. \\
$\hspace*{1em}$The regularity properties of the solutions to PDEs
as the one in (\ref{eq:elliptic}) have been widely investigated.
Actually, as we have already said, these PDEs can be considered as
the Euler-Lagrange equation of the functional appearing in (\ref{eq:P2})
\[
\mathcal{F}\left(u,\Omega\right)=\int_{\Omega}\frac{1}{q}\,(\vert\nabla u\vert-1)_{+}^{q}\,dx-\int_{\Omega}uf\,dx,
\]
which fits into the broader context of functionals having a \textit{$q$-Laplacian
type structure at infinity}. The local Lipschitz regularity for the
minimizers of asymptotically convex problems has been proved when
$f=0$ under quadratic growth assumptions in \cite{Chipot}, and later
in \cite{Giaq} for superquadratic growth conditions and in \cite{Leone}
in the case of subquadratic growth (see also \cite{Br00} and \cite{Fon}).
In addition, similar results have been obtained for functionals whose
integrands also depend on lower-order terms (see \cite{Ele2,Ele3,Pass}
and \cite{Ray} for example).\\
$\hspace*{1em}$A further topic concerning the regularity of solutions
to equations as the one in (\ref{eq:elliptic}) is the investigation
of their higher differentiability of both integer and fractional order,
and several results are available so far (see, for example, \cite{Bai,Clop}
and \cite{Gent}). In this regard, it must be mentioned that to the
best of our knowledge, different authors have proved the (local) higher
differentiability of integer order of a certain function of the gradient
of the solutions to (\ref{eq:elliptic}) that takes into account the
degeneracy of the equation in the case $q\geq2$ and $f\in W_{loc}^{1,p}\left(\Omega\right)$
: for further information see \cite{Br0,Br,Clop1}, as well as \cite{Cup}
for the case $f=0$. \\
$\hspace*{1em}$The technical difficulties encountered in dealing
with the singular case $1<q<2$ prevent us from achieving the same
kind of regularity results as those obtained for (\ref{eq:P1}) and
(\ref{eq:elliptic}) by Brasco, Carlier and Santambrogio in \cite{Br}:
there, they confine their analysis to the non-singular case $q\geq2$,
by assuming that $f\in W^{1,p}\left(\Omega\right)$. In particular,
they show that if $u\in W^{1,q}\left(\Omega\right)$ is a local weak
solution of (\ref{eq:elliptic}), then $\sigma_{0}=\nabla\mathcal{H}^{*}(\nabla u)\in W_{loc}^{1,r}\left(\Omega\right)$
for a suitable exponent $r=r(n,q)$. Moreover, using the Sobolev embedding
theorems, they find that $\sigma_{0}\in L_{loc}^{s}\left(\Omega\right)$
for all $s\in\left[1,+\infty\right)$ when $n=2$. As for the higher
integrability of $\nabla u$, they get $\nabla u\in L_{loc}^{q\,\frac{n}{n-2}}\left(\Omega\right)$
in the case $n>2$, while if $n=2$ they get $\nabla u\in L_{loc}^{s}\left(\Omega\right)$
for all $s\in\left[1,+\infty\right)$ (see \cite{Br}).\\
$\hspace*{1em}$However, despite the difficulties we have mentioned,
we are able to establish the following\\

\end{singlespace}
\begin{thm}
\begin{singlespace}
\noindent \label{thm:mainth} Let $n\geq2$, $q\in\left(1,2\right)$
and $\alpha\in\left(0,1\right)$. Moreover, let $u\in W_{loc}^{1,q}\left(\Omega\right)$
be a local weak solution of the equation
\begin{equation}
\mathrm{-div}\,(\nabla\mathcal{H}^{*}(\nabla u))=f.\label{eq:div}
\end{equation}
Then, for the function 
\[
H_{\frac{q}{2}}(\nabla u):=\left(\left|\nabla u\right|-1\right)_{+}^{\frac{q}{2}}\frac{\nabla u}{\left|\nabla u\right|}
\]
the following implications hold
\[
f\in W_{loc}^{1,p}\left(\Omega\right)\,\,\,\Rightarrow\,\,\,H_{\frac{q}{2}}(\nabla u)\in B_{2,\infty}^{\frac{1}{3-q}}\,\,\,\,locally\,\,in\,\,\Omega,
\]
\[
f\in B_{p,\infty,loc}^{\alpha}\left(\Omega\right)\,\,\,\Rightarrow\,\,\,H_{\frac{q}{2}}(\nabla u)\in B_{2,\infty}^{\min\left\{ \frac{\alpha+1}{2},\frac{1}{3-q}\right\} }\,\,\,\,locally\,\,in\,\,\Omega.
\]
Furthermore:\\
\\
$\mathrm{(}a\mathrm{)}$ if $f\in W_{loc}^{1,p}\left(\Omega\right)$,
then for any ball $B_{R}\Subset\Omega$ the following estimate 
\[
\underset{B_{R/2}}{\int}\left|\tau_{s,h}H_{\frac{q}{2}}(\nabla u)-H_{\frac{q}{2}}(\nabla u)\right|^{2}dx\leq C_{1}\,\left(\Vert\nabla f\Vert_{L^{p}\left(B'\right)}\,\Vert\nabla u\Vert_{L^{q}\left(B'\right)}\left|h\right|^{2}+\,\Vert\nabla u\Vert_{L^{q}\left(B'\right)}^{q}\left|h\right|^{\frac{2}{3-q}}\right)
\]
holds true for every $s\in\left\{ 1,\ldots,n\right\} $, for every
$h\in\mathbb{R}$ such that $\left|h\right|\leq r_{0}<\frac{1}{2}\,\mathrm{dist}\left(B_{R},\partial\Omega\right)$,
for $B'=B_{R}+B(0,r_{0})$ and a positive constant $C_{1}=C_{1}(R,q,n)$.\\

\noindent $\mathrm{(}b\mathrm{)}$ If, on the other hand, $f\in B_{p,\infty,loc}^{\alpha}\left(\Omega\right)$,
then for any ball $B_{R}\Subset\Omega$ the following estimate 
\[
\underset{B_{R/2}}{\int}\left|\tau_{s,h}H_{\frac{q}{2}}(\nabla u)-H_{\frac{q}{2}}(\nabla u)\right|^{2}dx\leq C_{2}\,\left(\Vert f\Vert_{B_{p,\infty}^{\alpha}\left(B'\right)}\,\Vert\nabla u\Vert_{L^{q}\left(B'\right)}\left|h\right|^{\alpha+1}+\,\Vert\nabla u\Vert_{L^{q}\left(B'\right)}^{q}\left|h\right|^{\frac{2}{3-q}}\right)
\]
holds true for every $s\in\left\{ 1,\ldots,n\right\} $, for every
$h\in\mathbb{R}$ such that $\left|h\right|\leq r_{0}<\frac{1}{2}\,\mathrm{dist}\left(B_{R},\partial\Omega\right)$,
for $B'=B_{R}+B(0,r_{0})$ and a positive constant $C_{2}=C_{2}(R,q,n)$. 
\end{singlespace}
\end{thm}

\begin{singlespace}
\noindent $\hspace*{1em}$The main tool in the proof of Theorem \ref{thm:mainth}
is the difference quotient technique, which, due to the singular character
of the case under consideration, requires new ideas with respect to
previous papers on the subject (\cite{Ace} and \cite{Tolk}), where
the case $1<q<2$ has been treated, but in the absence of singularity.
For the sake of clarity, let us now try to explain the difficulties
associated with the singular setting: when dealing with the higher
differentiability of solutions to subquadratic problems, the difference
quotient method usually yields a Caccioppoli-type inequality based
on the possibility of reabsorbing the terms in the right-hand side
that involve the increment of the gradient of a solution by the left-hand
side. In the singular setting, this possibility is actually excluded
by the fact that we cannot control the difference quotient of the
gradient of a solution $u$ in terms of $H_{\frac{q}{2}}(\nabla u)$,
whose difference quotient appears in the left-hand side, and this
is true regardless of the type of regularity of the datum $f$. Actually,
even in the case $f=0$, the integral appearing in the right-hand
side of formula \eqref{eq:stima0} cannot be estimated in a standard
way, since the difference quotient $\Delta_{h}u$ cannot be controlled
in terms of $\Delta_{h}(\nabla\mathcal{H}^{*}(\nabla u))$, unless
we perform the estimate far away from the set $\left\{ x:\vert\nabla u(x)\vert\leq1\right\} $.
On the other hand, as we have already observed, every $1$--Lipschitz
function $u$ is a solution of the homogeneous equation
\[
\mathrm{-div}\,(\nabla\mathcal{H}^{*}(\nabla u))=0,
\]
so that a more regular $f$ may not be sufficient to ensure a higher
differentiability property of the solutions.\\
$\hspace*{1em}$Returning to the variational problem (\ref{eq:P1}),
thanks to Theorem \ref{thm:mainth} we can obtain the following higher
differentiability result of fractional order for its (unique) minimizer:
\end{singlespace}
\begin{thm}
\begin{singlespace}
\noindent \label{thm:main2} Let $n\geq2$, $q\in\left(1,2\right)$
and $\alpha\in\left(0,1\right)$. Moreover, let $\sigma_{0}\in L^{p}\left(\Omega,\mathbb{R}^{n}\right)$
be the solution of $\left(\ref{eq:P1}\right)$. Then, the following
implications hold:
\[
f\in W_{loc}^{1,p}\left(\Omega\right)\,\,\,\Rightarrow\,\,\,\sigma_{0}\in B_{p,\infty}^{\frac{2}{p(3-q)}}\,\,\,\,locally\,\,in\,\,\Omega,
\]
\[
f\in B_{p,\infty,loc}^{\alpha}\left(\Omega\right)\,\,\,\Rightarrow\,\,\,\sigma_{0}\in B_{p,\infty}^{\min\left\{ \frac{\alpha+1}{p},\frac{2}{p(3-q)}\right\} }\,\,\,\,locally\,\,in\,\,\Omega.
\]
Furthermore:\\
\\
$\mathrm{(}a\mathrm{)}$ if $f\in W_{loc}^{1,p}\left(\Omega\right)$,
then for every solution $u$ of problem $\left(\ref{eq:P2}\right)$
and every ball $B_{R}\Subset\Omega$, the following estimate 
\[
\underset{B_{R/2}}{\int}\left|\tau_{s,h}\sigma_{0}-\sigma_{0}\right|^{p}dx\leq C_{1}\,\left(\Vert\nabla f\Vert_{L^{p}\left(B'\right)}\,\Vert\nabla u\Vert_{L^{q}\left(B'\right)}\left|h\right|^{2}+\,\Vert\nabla u\Vert_{L^{q}\left(B'\right)}^{q}\left|h\right|^{\frac{2}{3-q}}\right)
\]
holds true for every $s\in\left\{ 1,\ldots,n\right\} $, for every
$h\in\mathbb{R}$ such that $\left|h\right|\leq r_{0}<\frac{1}{2}\,\mathrm{dist}\left(B_{R},\partial\Omega\right)$,
for $B'=B_{R}+B(0,r_{0})$ and a positive constant $C_{1}=C_{1}(R,q,n)$.\\

\noindent $\mathrm{(}b\mathrm{)}$ If, on the other hand, $f\in B_{p,\infty,loc}^{\alpha}\left(\Omega\right)$,
then for every solution $u$ of problem $\left(\ref{eq:P2}\right)$
and every ball $B_{R}\Subset\Omega$, the following estimate 
\[
\underset{B_{R/2}}{\int}\left|\tau_{s,h}\sigma_{0}-\sigma_{0}\right|^{p}dx\leq C_{2}\,\left(\Vert f\Vert_{B_{p,\infty}^{\alpha}\left(B'\right)}\,\Vert\nabla u\Vert_{L^{q}\left(B'\right)}\left|h\right|^{\alpha+1}+\,\Vert\nabla u\Vert_{L^{q}\left(B'\right)}^{q}\left|h\right|^{\frac{2}{3-q}}\right)
\]
holds true for every $s\in\left\{ 1,\ldots,n\right\} $, for every
$h\in\mathbb{R}$ such that $\left|h\right|\leq r_{0}<\frac{1}{2}\,\mathrm{dist}\left(B_{R},\partial\Omega\right)$,
for $B'=B_{R}+B(0,r_{0})$ and a positive constant $C_{2}=C_{2}(R,q,n)$.
\end{singlespace}
\end{thm}

\begin{singlespace}
\noindent \begin{brem} We immediately observe that if we take $q=p=2$
and $f\in W_{loc}^{1,2}\left(\Omega\right)$ in the statements $\mathrm{(}a\mathrm{)}$
of Theorems \ref{thm:mainth} and \ref{thm:main2}, then we get back
the full Sobolev regularity results of \cite[Theorem 4.2 and Corollary 4.3]{Br},
i.e. there is no discrepancy between the old and the new results.
\end{brem}\bigskip{}

\noindent $\hspace*{1em}$As a consequence of the above theorem, using
a local version of Sobolev's embedding theorem for Besov spaces, we
obtain the following higher integrability result for $\sigma_{0}$
:
\end{singlespace}
\begin{cor}
\begin{singlespace}
\noindent \label{cor:corol1} Under the assumptions of Theorem \ref{thm:main2},
we obtain the following implications:
\[
f\in W_{loc}^{1,p}\left(\Omega\right)\,\,\,\Rightarrow\,\,\,\sigma_{0}\in L_{loc}^{r}\left(\Omega\right)\,\,\,for\,\,all\,\,r\in\left[1,\frac{np(3-q)}{n(3-q)-2}\right),
\]
\[
f\in B_{p,\infty,loc}^{\alpha}\left(\Omega\right)\,\,\,\Rightarrow\,\,\,\sigma_{0}\in L_{loc}^{s}\left(\Omega\right),
\]
where
\[
s=\begin{cases}
\begin{array}{cc}
any\,\,value<\frac{np}{n-\alpha-1},\,\, & if\,\,\,0<\alpha<\frac{q-1}{3-q},\\
\\
any\,\,value<\frac{np(3-q)}{n(3-q)-2}, & if\,\,\,\frac{q-1}{3-q}\leq\alpha<1.
\end{array}\end{cases}
\]
\end{singlespace}
\end{cor}

\begin{singlespace}
\noindent $\hspace*{1em}$Finally, as an easy consequence of the previous
results, we also get a gain of integrability for the gradient of the
local weak solutions of (\ref{eq:div}). Indeed, we have the following
\end{singlespace}
\begin{cor}
\begin{singlespace}
\noindent \label{cor:corol3} Under the assumptions of Theorem \ref{thm:mainth},
we obtain the following implications: 
\[
f\in W_{loc}^{1,p}\left(\Omega\right)\,\,\,\Rightarrow\,\,\,\nabla u\in L_{loc}^{r}\left(\Omega\right)\,\,\,for\,\,all\,\,r\in\left[1,\frac{nq(3-q)}{n(3-q)-2}\right),
\]
\[
f\in B_{p,\infty,loc}^{\alpha}\left(\Omega\right)\,\,\,\Rightarrow\,\,\,\nabla u\in L_{loc}^{s}\left(\Omega\right),
\]
where 
\[
s=\begin{cases}
\begin{array}{cc}
any\,\,value<\frac{nq}{n-\alpha-1},\,\, & if\,\,\,0<\alpha<\frac{q-1}{3-q},\\
\\
any\,\,value<\frac{nq(3-q)}{n(3-q)-2}, & if\,\,\,\frac{q-1}{3-q}\leq\alpha<1.
\end{array}\end{cases}
\]
\end{singlespace}
\end{cor}

\begin{singlespace}
\noindent $\hspace*{1em}$For the sake of completeness, in this paper
we shall prove that similar results also hold in the case $q\geq2$,
by assuming that $f\in B_{p,\infty,loc}^{\alpha}\left(\Omega\right)$
for some $\alpha\in\left(0,1\right)$. Actually, on the one hand,
Theorem \ref{thm:q>2} below extends the results proved in \cite{Clop}
to a widely degenerate setting; on the other hand, it extends the
aforementioned results contained in \cite{Br} to the case of data
in a local Besov space.\\
$\hspace*{1em}$It is worth mentioning that if $f\in W_{\diamond}^{1,p}\left(\Omega\right)$,
then for suitable values of $n\geq2$ and $p>2$ there are some cases
in which it is possible to prove that every solution $u\in W_{\diamond}^{1,q}\left(\Omega\right)$
of the Neumann boundary problem (\ref{eq:Eul-Lag}) belongs to $W^{1,\infty}\left(\Omega\right)$
(see \cite[Theorem 2.1]{Br00}). However, in these cases, the boundedness
of $\nabla u$ gives neither a weak differentiability of integer order
for $H_{\frac{q}{2}}(\nabla u)$ nor a fractional differentiability
of higher order than that obtained in Theorem \ref{thm:mainth}.\\
$\hspace*{1em}$Before describing the structure of this paper, we
wish to briefly motivate what follows by pointing out that problems
such as (\ref{eq:P1}) and (\ref{eq:P2}) and degenerate equations
of the form (\ref{eq:elliptic}) may arise in \textit{optimal transport
problems with congestion effects} on very dense two-dimensional networks.
These problems lead, at the limit, to minimization problems posed
on measures on curves, as shown in \cite{Baillon}. In recent works
(\cite{Br0} and \cite{Br} again), Brasco, Carlier and Santambrogio
have gone one step further, by showing that these latter problems
can be reformulated in terms of the minimization of an integral functional
over a set of vector fields with prescribed divergence, exactly as
in problem (\ref{eq:P1}), which was introduced by Beckmann in the
1950s (see \cite{Beck}). For further details, we refer the interested
reader to \cite[Sections 1-2]{Br00} and \cite[Section 1]{Br1}.\\
\\
$\hspace*{1em}$The paper is organized as follows. Section \ref{sec:prelim}
is devoted to the preliminaries: after a list of some classical notations
and some essential estimates, we recall the basic properties of Besov-Lipschitz
spaces and some elementary properties of the difference quotients
of Sobolev functions. In Section \ref{sec:1<q<2}, we prove Theorems
\ref{thm:mainth} and \ref{thm:main2} and Corollaries \ref{cor:corol1}
and \ref{cor:corol3}. Finally, in Section \ref{sec:q>2} we give
similar results in the non-singular case $q\geq2$.
\end{singlespace}
\begin{singlespace}

\section{Notations and preliminaries\label{sec:prelim}}
\end{singlespace}

\begin{singlespace}
\noindent $\hspace*{1em}$In this paper we shall denote by $C$ or
$c$ a general positive constant that may vary on different occasions.
Relevant dependencies on parameters and special constants will be
suitably emphasized using parentheses or subscripts. The norm we use
on $\mathbb{R}^{n}$ will be the standard Euclidean one and it will
be denoted by $\left|\,\cdot\,\right|$. In particular, for the vectors
$\xi,\eta\in\mathbb{R}^{n}$, we write $\langle\xi,\eta\rangle$ for
the usual inner product and $\left|\xi\right|:=\langle\xi,\xi\rangle^{\frac{1}{2}}$
for the corresponding Euclidean norm.\\
$\hspace*{1em}$In what follows, $B(x,r)=B_{r}(x)=\left\{ y\in\mathbb{R}^{n}:\left|y-x\right|<r\right\} $
will denote the ball centered at $x$ of radius $r$. We shall omit
the dependence on the center and on the radius when no confusion arises.
\\
$\hspace*{1em}$We now recall some tools that will be useful to prove
our results. For the auxiliary function $V_{\gamma}$, defined for
all $\xi\in\mathbb{R}^{n}$ as
\[
V_{\gamma}(\xi):=\left(\mu^{2}+\left|\xi\right|^{2}\right)^{\gamma/2}\xi,
\]
where $\mu\geq0$ and $\gamma>0$ are parameters, we record the following
estimates (see \cite{Giaq}):\\

\end{singlespace}
\begin{lem}
\begin{singlespace}
\noindent \label{lem:Giaquinta} Let $\gamma>0$ and $\mu\geq0$.
Then there exist two positive constants $c_{1}(\gamma)$ and $c_{2}(\gamma)$
such that 
\[
c_{1}(\gamma)\,\left|\xi-\eta\right|\leq\frac{\left|V_{\gamma}(\xi)-V_{\gamma}(\eta)\right|}{\left(\mu^{2}+\left|\xi\right|^{2}+\left|\eta\right|^{2}\right)^{\gamma/2}}\leq c_{2}(\gamma)\,\left|\xi-\eta\right|
\]
for any $\xi,\eta\in\mathbb{R}^{n}$. 
\end{singlespace}
\end{lem}

\begin{singlespace}
\noindent $\hspace*{1em}$The next result is an extension of the previous
lemma to the case $\gamma\in\left(-1,0\right)$ (see the proof of
\cite[Lemma 2.2]{Ace}):
\end{singlespace}
\begin{lem}
\begin{singlespace}
\noindent \label{lem:Fusco} For every $\gamma\in\left(-1/2,0\right)$
and $\mu\geq0$ we have 
\begin{equation}
\left(2\gamma+1\right)\left|\xi-\eta\right|\leq\frac{\left|\left(\mu^{2}+\left|\xi\right|^{2}\right)^{\gamma}\xi-\left(\mu^{2}+\left|\eta\right|^{2}\right)^{\gamma}\eta\right|}{\left(\mu^{2}+\left|\xi\right|^{2}+\left|\eta\right|^{2}\right)^{\gamma}}\leq\frac{c(n)}{2\gamma+1}\left|\xi-\eta\right|\label{eq:Acerbi}
\end{equation}
for any $\xi,\eta\in\mathbb{R}^{n}$. Moreover, setting
\[
F(\zeta):=\frac{1}{2(\gamma+1)}\left(\mu^{2}+\left|\zeta\right|^{2}\right)^{\gamma+1},\,\,\,\,\,\zeta\in\mathbb{R}^{n},
\]

\noindent we obtain
\begin{equation}
\nabla F(\zeta)=\left(\mu^{2}+\left|\zeta\right|^{2}\right)^{\gamma}\zeta,\,\,\,\,\,\zeta\in\mathbb{R}^{n},\label{eq:Ace1}
\end{equation}
and
\begin{equation}
\langle\nabla F(\xi)-\nabla F(\eta),\xi-\eta\rangle\geq(2\gamma+1)\left(\mu^{2}+\left|\xi\right|^{2}+\left|\eta\right|^{2}\right)^{\gamma}\left|\xi-\eta\right|^{2}\label{eq:Ace2}
\end{equation}
for any $\xi,\eta\in\mathbb{R}^{n}$.
\end{singlespace}
\end{lem}

\begin{singlespace}
\noindent $\hspace*{1em}$In the following, we shall also use the
auxiliary function $H_{\alpha}:\mathbb{R}^{n}\rightarrow\mathbb{R}^{n}$
defined as
\[
H_{\alpha}(\xi):=\left(\left|\xi\right|-1\right)_{+}^{\alpha}\frac{\xi}{\left|\xi\right|},
\]
where $\alpha>0$ is a parameter. As for the above function, we can
obtain the following result that will be needed to prove Theorem \ref{thm:mainth}
:
\end{singlespace}
\begin{lem}
\begin{singlespace}
\noindent \label{lem:jolly} For every $\left(\alpha,\varepsilon\right)\in\mathbb{R}^{+}\times\mathbb{R}^{+}$
with $\alpha<\varepsilon$, there exist two positive constants $\beta_{1}(\alpha,\varepsilon)$
and $\beta_{2}(\alpha,\varepsilon,n)$ such that
\[
\beta_{1}\,\left|H_{\varepsilon}(\xi)-H_{\varepsilon}(\eta)\right|\leq\frac{\left|H_{\alpha}(\xi)-H_{\alpha}(\eta)\right|}{\left(\left(\left|\xi\right|-1\right)_{+}^{\varepsilon}+\left(\left|\eta\right|-1\right)_{+}^{\varepsilon}\right)^{\frac{\alpha-\varepsilon}{\varepsilon}}}\leq\beta_{2}\,\left|H_{\varepsilon}(\xi)-H_{\varepsilon}(\eta)\right|
\]
for any $\xi,\eta\in\mathbb{R}^{n}$.
\end{singlespace}
\end{lem}

\begin{singlespace}
\noindent \begin{proof}[\bfseries{Proof}] By Lemma \ref{lem:Fusco},
we have that for every $\gamma\in\left(-1,0\right)$ there exist two
positive constants $k_{1}(\gamma)$ and $k_{2}(\gamma,n)$ such that
\begin{equation}
k_{1}\,\left|z-w\right|^{2}\leq\frac{\left|\left|z\right|^{\gamma}z-\left|w\right|^{\gamma}w\right|^{2}}{\left(\left|z\right|^{2}+\left|w\right|^{2}\right)^{\gamma}}\leq k_{2}\,\left|z-w\right|^{2}\label{eq:GiaFus}
\end{equation}
for any $z,w\in\mathbb{R}^{n}$. Now we recall that for every $\gamma\in\left(-1,0\right)$
and every $z,w\in\mathbb{R}^{n}$ we get
\[
\frac{1}{\left(\left|z\right|^{2}+\left|w\right|^{2}\right)^{\gamma}}\leq\frac{1}{(\left|z\right|+\left|w\right|)^{2\gamma}}\,\,\,\,\,\,\,\,\mathrm{and}\,\,\,\,\,\,\,\,\frac{1}{\left(\left|z\right|^{2}+\left|w\right|^{2}\right)^{\gamma}}\geq\frac{2^{\gamma}}{(\left|z\right|+\left|w\right|)^{2\gamma}}.
\]
From this observation and from (\ref{eq:GiaFus}), it follows that
for every $\gamma\in\left(-1,0\right)$ there exist two positive constants
$k_{3}(\gamma)$ and $k_{4}(\gamma,n)$ such that
\begin{equation}
k_{3}\,\left|z-w\right|\leq\frac{\left|\left|z\right|^{\gamma}z-\left|w\right|^{\gamma}w\right|}{(\left|z\right|+\left|w\right|)^{\gamma}}\leq k_{4}\,\left|z-w\right|\label{eq:GiFu2}
\end{equation}
for any $z,w\in\mathbb{R}^{n}$. Then, for any fixed $\left(\alpha,\varepsilon\right)\in\mathbb{R}^{+}\times\mathbb{R}^{+}$
such that $\alpha<\varepsilon$, by taking
\[
z=H_{\varepsilon}(\xi),\,\,\,\,\,\,w=H_{\varepsilon}(\eta)\,\,\,\,\,\,\mathrm{and}\,\,\,\,\,\,\gamma=\frac{\alpha-\varepsilon}{\varepsilon}
\]
into (\ref{eq:GiFu2}), we obtain the desired conclusion.\end{proof}

\noindent $\hspace*{1em}$Now we recall a result that will be needed
to prove Theorem \ref{thm:q>2} and whose proof can be found in \cite[Lemma 4.1]{Br}
:
\end{singlespace}
\begin{lem}
\begin{singlespace}
\noindent \label{lem:Brasco} If $2\leq q<\infty$, then for every
$\xi,\eta\in\mathbb{R}^{n}$ we get 
\[
\langle H_{q-1}(\xi)-H_{q-1}(\eta),\xi-\eta\rangle\geq\frac{4}{q^{2}}\left|H_{\frac{q}{2}}(\xi)-H_{\frac{q}{2}}(\eta)\right|^{2},
\]
\[
\left|H_{q-1}(\xi)-H_{q-1}(\eta)\right|\leq(q-1)\left(\left|H_{\frac{q}{2}}(\xi)\right|^{\frac{q-2}{q}}+\left|H_{\frac{q}{2}}(\eta)\right|^{\frac{q-2}{q}}\right)\left|H_{\frac{q}{2}}(\xi)-H_{\frac{q}{2}}(\eta)\right|.
\]
\end{singlespace}
\end{lem}

\begin{singlespace}
\noindent $\hspace*{1em}$For further needs, we shall now prove that
the first estimate of the previous lemma holds true also in the case
$1<q<2$. More precisely, we have the following 
\end{singlespace}
\begin{lem}
\begin{singlespace}
\noindent \label{lem:below} If $1<q<2$, then there exists a constant
$\beta\equiv\beta(q,n)>0$ such that 
\[
\langle H_{q-1}(\xi)-H_{q-1}(\eta),\xi-\eta\rangle\geq\beta\left|H_{\frac{q}{2}}(\xi)-H_{\frac{q}{2}}(\eta)\right|^{2}
\]
for every $\xi,\eta\in\mathbb{R}^{n}$. 
\end{singlespace}
\end{lem}

\begin{singlespace}
\noindent \begin{proof}[\bfseries{Proof}] We first note that the
above inequality is trivially satisfied when $\left|\xi\right|,\left|\eta\right|\leq1$.
Taking $\gamma=(q-2)/2$ and $\mu=0$, equality (\ref{eq:Ace1}) becomes
\[
\nabla F(\zeta)=\left|\zeta\right|^{q-2}\zeta,\,\,\,\,\,\zeta\in\mathbb{R}^{n},
\]
and estimates (\ref{eq:Ace2}) and (\ref{eq:Acerbi}) imply \begin{equation}\label{eq:ACEFUSCO}
\begin{split}
\langle\left|z\right|^{q-2}z-\left|w\right|^{q-2}w,z-w\rangle
&\geq(q-1)\left(\left|z\right|^{2}+\left|w\right|^{2}\right)^{\frac{q-2}{2}}\left|z-w\right|^{2}\\
&\geq\beta(q,n)\left|\left|z\right|^{\frac{q-2}{2}}z-\left|w\right|^{\frac{q-2}{2}}w\right|^{2}
\end{split}
\end{equation}for any $z,w\in\mathbb{R}^{n}$. Taking $z=(\left|\xi\right|-1)_{+}\,\xi/\left|\xi\right|$
and $w=(\left|\eta\right|-1)_{+}\,\eta/\left|\eta\right|$ into \eqref{eq:ACEFUSCO},
we obtain 
\begin{equation}
\bigg\langle H_{q-1}(\xi)-H_{q-1}(\eta),(\left|\xi\right|-1)_{+}\,\frac{\xi}{\left|\xi\right|}-(\left|\eta\right|-1)_{+}\,\frac{\eta}{\left|\eta\right|}\bigg\rangle\geq\beta(q,n)\left|H_{\frac{q}{2}}(\xi)-H_{\frac{q}{2}}(\eta)\right|^{2}.\label{eq:disFusco4}
\end{equation}
$\hspace*{1em}$Now let $\xi,\eta\in\mathbb{R}^{n}$ be such that
$\left|\xi\right|,\left|\eta\right|>1$ and write the left-hand side
of the previous inequality as the difference of two terms $A-B$,
where
\[
A:=\langle H_{q-1}(\xi)-H_{q-1}(\eta),\xi-\eta\rangle
\]
and
\[
B:=\bigg\langle H_{q-1}(\xi)-H_{q-1}(\eta),\xi-(\left|\xi\right|-1)_{+}\,\frac{\xi}{\left|\xi\right|}-\eta+(\left|\eta\right|-1)_{+}\,\frac{\eta}{\left|\eta\right|}\bigg\rangle.
\]
\\
Using the fact that $\left|\xi\right|,\left|\eta\right|>1$ and the
Cauchy-Schwarz inequality, we get \begin{align*}
B&=\bigg\langle H_{q-1}(\xi)-H_{q-1}(\eta),\frac{\xi}{\left|\xi\right|}-\frac{\eta}{\left|\eta\right|}\bigg\rangle\\
&=(\left|\xi\right|-1)^{q-1}+(\left|\eta\right|-1)^{q-1}-\frac{\langle\xi,\eta\rangle}{\left|\xi\right|\left|\eta\right|}\left[(\left|\xi\right|-1)^{q-1}+(\left|\eta\right|-1)^{q-1}\right]\\
&\geq(\left|\xi\right|-1)^{q-1}+(\left|\eta\right|-1)^{q-1}-\frac{\left|\xi\right|\left|\eta\right|}{\left|\xi\right|\left|\eta\right|}\left[(\left|\xi\right|-1)^{q-1}+(\left|\eta\right|-1)^{q-1}\right]=0,
\end{align*}from which the assertion directly follows. \\
$\hspace*{1em}$Finally, let $\left|\xi\right|>1$ and $\left|\eta\right|\leq1$.
In this case, we have $H_{q-1}(\eta)=H_{\frac{q}{2}}(\eta)=0$ and
we can write the left-hand side of (\ref{eq:disFusco4}) as the difference
of two terms $D-E$, where
\[
D:=\langle H_{q-1}(\xi),\xi-\eta\rangle\,\,\,\,\,\,\,\,\mathrm{and}\,\,\,\,\,\,\,\,E:=\bigg\langle H_{q-1}(\xi),\frac{\xi}{\left|\xi\right|}-\eta\bigg\rangle.
\]
Using the fact that $\left|\xi\right|>1,\left|\eta\right|\leq1$ and
the Cauchy-Schwarz inequality, we find that
\[
E=(\left|\xi\right|-1)^{q-1}-\frac{\langle\xi,\eta\rangle}{\left|\xi\right|}(\left|\xi\right|-1)^{q-1}\geq(\left|\xi\right|-1)^{q-1}-\frac{\left|\xi\right|\left|\eta\right|}{\left|\xi\right|}(\left|\xi\right|-1)^{q-1}\geq0,
\]
which concludes the proof.\end{proof} 
\end{singlespace}
\begin{singlespace}

\subsection{Besov--Lipschitz spaces}
\end{singlespace}

\begin{singlespace}
\noindent $\hspace*{1em}$Given $h\in\mathbb{R}^{n}$ and $v:\mathbb{R}^{n}\rightarrow\mathbb{R}$,
let us introduce the notations $\tau_{h}v(x)=v(x+h)$ and $\Delta_{h}[v](x)=v(x+h)-v(x)$.
As in \cite[Section 2.5.12]{Tri}, given $0<\alpha<1$ and $1\leq p,q<\infty$,
we say that $v$ belongs to the Besov space $B_{p,q}^{\alpha}\left(\mathbb{R}^{n}\right)$
if $v\in L^{p}\left(\mathbb{R}^{n}\right)$ and
\begin{equation}
\left[v\right]_{\dot{B}_{p,q}^{\alpha}\left(\mathbb{R}^{n}\right)}:=\left(\int_{\mathbb{R}^{n}}\left(\int_{\mathbb{R}^{n}}\frac{\left|\Delta_{h}[v](x)\right|^{p}}{\left|h\right|^{\alpha p}}\,dx\right)^{\frac{q}{p}}\frac{dh}{\left|h\right|^{n}}\right)^{\frac{1}{q}}<\infty.\label{eq:Bes1}
\end{equation}
One can define a norm on the space $B_{p,q}^{\alpha}\left(\mathbb{R}^{n}\right)$
as follows
\[
\Vert v\Vert_{B_{p,q}^{\alpha}\left(\mathbb{R}^{n}\right)}:=\Vert v\Vert_{L^{p}\left(\mathbb{R}^{n}\right)}+\left[v\right]_{\dot{B}_{p,q}^{\alpha}\left(\mathbb{R}^{n}\right)},
\]
and with this norm $B_{p,q}^{\alpha}\left(\mathbb{R}^{n}\right)$
is a Banach space. Equivalently, we could simply say that a function
$v\in L^{p}\left(\mathbb{R}^{n}\right)$ belongs to $B_{p,q}^{\alpha}\left(\mathbb{R}^{n}\right)$
if and only if $\frac{\Delta_{h}[v]}{\left|h\right|^{\alpha}}\in L^{q}\left(\frac{dh}{\left|h\right|^{n}};L^{p}\left(\mathbb{R}^{n}\right)\right)$.
As usual, in (\ref{eq:Bes1}) if one simply integrates for $h\in B(0,\delta)$
for a fixed $\delta>0$, then an equivalent norm is obtained, since
\[
\left(\int_{\left\{ \left|h\right|\geq\delta\right\} }\left(\int_{\mathbb{R}^{n}}\frac{\left|\Delta_{h}[v](x)\right|^{p}}{\left|h\right|^{\alpha p}}\,dx\right)^{\frac{q}{p}}\frac{dh}{\left|h\right|^{n}}\right)^{\frac{1}{q}}\leq c(n,\alpha,p,q,\delta)\,\Vert v\Vert_{L^{p}\left(\mathbb{R}^{n}\right)}.
\]
Similarly, for a function $v\in L^{p}\left(\mathbb{R}^{n}\right)$
we say that $v\in B_{p,\infty}^{\alpha}\left(\mathbb{R}^{n}\right)$
if 
\begin{equation}
\left[v\right]_{\dot{B}_{p,\infty}^{\alpha}\left(\mathbb{R}^{n}\right)}:=\sup_{h\in\mathbb{R}^{n}}\left(\int_{\mathbb{R}^{n}}\frac{\left|\Delta_{h}[v](x)\right|^{p}}{\left|h\right|^{\alpha p}}\,dx\right)^{\frac{1}{p}}<\infty,\label{eq:Bes2}
\end{equation}
and we can define the following norm
\[
\Vert v\Vert_{B_{p,\infty}^{\alpha}\left(\mathbb{R}^{n}\right)}:=\Vert v\Vert_{L^{p}\left(\mathbb{R}^{n}\right)}+\left[v\right]_{\dot{B}_{p,\infty}^{\alpha}\left(\mathbb{R}^{n}\right)}.
\]
Again, in (\ref{eq:Bes2}) one can simply take the supremum over $\left|h\right|\leq\delta$
for a fixed $\delta>0$, thus obtaining an equivalent norm. By construction,
$B_{p,q}^{\alpha}\left(\mathbb{R}^{n}\right)\subset L^{p}\left(\mathbb{R}^{n}\right)$.
Moreover, one also has the following version of Sobolev embeddings
(a proof can be found in \cite[Proposition 7.12]{Har}, taking into
account that $L^{r}=F_{r,2}^{0}$, with $1<r<+\infty$):\\

\end{singlespace}
\begin{lem}
\begin{singlespace}
\noindent \label{lem:BesEmbed} Suppose that $0<\alpha<1$.\\
 \\
$\mathrm{(}a\mathrm{)}$ If $1<p<\frac{n}{\alpha}$ and $1\leq q\leq p_{\alpha}^{*}:=\frac{np}{n-\alpha p}$,
then there exists a continuous embedding $B_{p,q}^{\alpha}\left(\mathbb{R}^{n}\right)\hookrightarrow$\\
\textcolor{white}{$\mathrm{(a)}\,$}$L^{p_{\alpha}^{*}}\left(\mathbb{R}^{n}\right)$.\\
$\mathrm{(}b\mathrm{)}$ If $p=\frac{n}{\alpha}$ and $1\leq q\leq\infty$,
then there exists a continuous embedding $B_{p,q}^{\alpha}\left(\mathbb{R}^{n}\right)\hookrightarrow BMO\left(\mathbb{R}^{n}\right)$,\\
\\
where $BMO$ denotes the space of functions with bounded mean oscillations
\cite[Chapter 2]{Giu}.
\end{singlespace}
\end{lem}

\begin{singlespace}
\noindent $\hspace*{1em}$For further needs, we recall the following
inclusions (see \cite[Proposition 7.10 and Formula (7.35)]{Har}).
\end{singlespace}
\begin{lem}
\begin{singlespace}
\noindent \label{lem:BesInclu} Suppose that $0<\beta<\alpha<1$.\\
\\
$\mathrm{(}a\mathrm{)}$ If $1<p\leq+\infty$ and $1\leq q\leq r\leq+\infty$
then $B_{p,q}^{\alpha}\left(\mathbb{R}^{n}\right)\subset B_{p,r}^{\alpha}\left(\mathbb{R}^{n}\right)$.\\
$\mathrm{(}b\mathrm{)}$ If $1<p\leq+\infty$ and $1\leq q$, $r\leq+\infty$
then $B_{p,q}^{\alpha}\left(\mathbb{R}^{n}\right)\subset B_{p,r}^{\beta}\left(\mathbb{R}^{n}\right)$.\\
$\mathrm{(}c\mathrm{)}$ If $1\leq q\leq+\infty$, then $B_{\frac{n}{\alpha},q}^{\alpha}\left(\mathbb{R}^{n}\right)\subset B_{\frac{n}{\beta},q}^{\beta}\left(\mathbb{R}^{n}\right)$.
\end{singlespace}
\end{lem}

\begin{singlespace}
\noindent $\hspace*{1em}$Combining Lemmas \ref{lem:BesEmbed} and
\ref{lem:BesInclu}, we obtain the following Sobolev-type embedding
theorem for Besov spaces $B_{p,\infty}^{\alpha}\left(\mathbb{R}^{n}\right)$
that are excluded from assumptions $\mathrm{(}a\mathrm{)}$ in Lemma
\ref{lem:BesEmbed}.
\end{singlespace}
\begin{thm}
\begin{singlespace}
\noindent Suppose that $0<\beta<\alpha<1$. If $1<p<\frac{n}{\alpha}$,
then there exists a continuous embedding $B_{p,\infty}^{\alpha}\left(\mathbb{R}^{n}\right)\hookrightarrow L^{p_{\beta}^{*}}\left(\mathbb{R}^{n}\right)$.
Moreover, for every $v\in B_{p,\infty}^{\alpha}\left(\mathbb{R}^{n}\right)$,
the following local estimate 
\[
\Vert v\Vert_{L^{\frac{np}{n-\beta p}}\left(B_{\rho}\right)}\leq c\left(\left[v\right]_{\dot{B}_{p,\infty}^{\alpha}\left(B_{R}\right)}+\Vert v\Vert_{L^{p}\left(B_{R}\right)}\right)
\]
holds for every ball $B_{\rho}\subset B_{R}$ with $c=c\,(n,R,\rho,\alpha,\beta)$.
\end{singlespace}
\end{thm}

\begin{singlespace}
\noindent $\hspace*{1em}$We can also define local Besov spaces as
follows. Given a domain $\Omega\subset\mathbb{R}^{n}$, we say that
a function $v$ belongs to $B_{p,q,loc}^{\alpha}\left(\Omega\right)$
if $\varphi v\in B_{p,q}^{\alpha}\left(\mathbb{R}^{n}\right)$ whenever
$\varphi$ belongs to the class $C_{c}^{\infty}\left(\Omega\right)$
of smooth functions with compact support contained in $\Omega$. It
is worth noticing that one can prove suitable versions of Lemmas \ref{lem:BesEmbed}
and \ref{lem:BesInclu}, by using local Besov spaces.\\
$\hspace*{1em}$The following lemma is an easy consequence of the
definitions given above and its proof can be found in \cite{Bai}. 
\end{singlespace}
\begin{lem}
\begin{singlespace}
\noindent A function $v\in L_{loc}^{p}\left(\Omega\right)$ belongs
to the local Besov space $B_{p,q,loc}^{\alpha}\left(\Omega\right)$
if and only if 
\[
\bigg\Vert\frac{\Delta_{h}[v]}{\left|h\right|^{\alpha}}\bigg\Vert_{L^{q}\left(\frac{dh}{\left|h\right|^{n}};L^{p}\left(B\right)\right)}<\infty
\]
for any ball $B\subset B_{2\,r_{B}}\Subset\Omega$ with radius $r_{B}$.
Here, the balls $B$ and $B_{2\,r_{B}}$ are supposed to be concentric
and the measure $\frac{dh}{\left|h\right|^{n}}$ is restricted to
the ball $B(0,r_{B})$ on the $h$-space.
\end{singlespace}
\end{lem}

\begin{singlespace}
\noindent $\hspace*{1em}$Finally, the following result represents,
in some sense, the local counterpart of Lemma \ref{lem:BesEmbed}
in the case $q=\infty$ :
\end{singlespace}
\begin{thm}
\begin{singlespace}
\noindent \label{thm:emb2} On any domain $\Omega\subset\mathbb{R}^{n}$
we have the continuous embedding $B_{p,\infty,loc}^{\alpha}\left(\Omega\right)\hookrightarrow L_{loc}^{r}\left(\Omega\right)$
for all $r<\frac{np}{n-\alpha p}$, provided $\alpha\in\left(0,1\right)$
and $1<p<\frac{n}{\alpha}$. 
\end{singlespace}
\end{thm}

\begin{singlespace}
\noindent We refer to \cite[Sections 30-32]{Tar} for a proof of this
theorem. In fact, the above statement follows by localizing the corresponding
result proved for functions defined on $\mathbb{R}^{n}$ in \cite{Tar},
by simply using a smooth cut off function. 
\end{singlespace}
\begin{singlespace}

\subsection{Difference quotients }
\end{singlespace}

\begin{singlespace}
\noindent $\hspace*{1em}$We recall here the definition and some elementary
properties of the difference quotients that will be useful in the
following (see, for example, \cite{Giu}).\\

\end{singlespace}
\begin{defn}
\begin{singlespace}
\noindent Given $h\in\mathbb{R}\setminus\left\{ 0\right\} $, for
every vector-valued function $F:\Omega\rightarrow\mathbb{R}^{N}$
defined in an open set $\Omega\subset\mathbb{R}^{n}$, we call the
\textit{difference quotient} of $F$ with respect to $x_{s}$ the
function
\[
\Delta_{s,h}F(x)=\frac{\tau_{s,h}F(x)-F(x)}{h}=\frac{F(x+he_{s})-F(x)}{h},
\]
where $e_{s}$ is the unit vector in the $x_{s}$ direction and $s\in\left\{ 1,\ldots,n\right\} $.
\end{singlespace}
\end{defn}

\begin{singlespace}
\noindent $\hspace*{1em}$When no confusion can arise, we shall omit
the index $s$ and simply write $\Delta_{h}$ or $\tau_{h}$ instead
of $\Delta_{s,h}$ or $\tau_{s,h}$, respectively.
\end{singlespace}
\begin{prop}
\begin{singlespace}
\noindent Let $F$ and $G$ be two functions such that $F,G\in W^{1,p}\left(\Omega\right)$,
with $p\geq1$, and let us consider the set
\[
\Omega_{\left|h\right|}:=\left\{ x\in\Omega:\mathrm{dist}\left(x,\partial\Omega\right)>\left|h\right|\right\} .
\]
Then:\\
\\
$\mathrm{(}a\mathrm{)}$ $\Delta_{h}F\in W^{1,p}\left(\Omega_{\left|h\right|}\right)$
and
\[
D_{i}(\Delta_{h}F)=\Delta_{h}(D_{i}F),\,\,\,\,\,for\,\,every\,\,i\in\left\{ 1,\ldots,n\right\} .
\]
$\mathrm{(}b\mathrm{)}$ If at least one of the functions $F$ or
$G$ has support contained in $\Omega_{\left|h\right|}$, then
\[
\int_{\Omega}F\,\Delta_{h}G\,dx=-\int_{\Omega}G\,\Delta_{-h}F\,dx.
\]
$\mathrm{(}c\mathrm{)}$ We have 
\[
\Delta_{h}(FG)(x)=F(x+he_{s})\Delta_{h}G(x)\,+\,G(x)\Delta_{h}F(x).
\]
\end{singlespace}
\end{prop}

\begin{singlespace}
\noindent $\hspace*{1em}$The next result is a kind of integral version
of Lagrange Theorem and it will be used to prove Theorems \ref{thm:mainth}
and \ref{thm:q>2}.
\end{singlespace}
\begin{lem}
\begin{singlespace}
\noindent \label{lem:diffquo} If $0<\rho<R$, $\left|h\right|<\frac{R-\rho}{2}$,
$1<p<+\infty$, and $F\in L^{p}\left(B_{R},\mathbb{R}^{N}\right)$,
$DF\in L^{p}\left(B_{R},\mathbb{R}^{N\times n}\right)$, then
\[
\int_{B_{\rho}}\left|F(x+he_{s})-F(x)\right|^{p}dx\leq c^{p}(n)\left|h\right|^{p}\int_{B_{R}}\left|DF(x)\right|^{p}dx,
\]
where the balls $B_{\rho}$ and $B_{R}$ are supposed to be concentric.
Moreover
\[
\int_{B_{\rho}}\left|F(x+he_{s})\right|^{p}dx\leq\int_{B_{R}}\left|F(x)\right|^{p}dx.
\]
\end{singlespace}
\end{lem}

\begin{singlespace}

\section{The singular case $1<q<2$\label{sec:1<q<2}}
\end{singlespace}

\begin{singlespace}
\noindent $\hspace*{1em}$Having declared our aims and introduced
all the required tools, we now come to the\\

\noindent \begin{proof}[\bfseries{Proof of Theorem~\ref{thm:mainth}}]
We first observe that $H_{\frac{q}{2}}(\nabla u)\in L_{loc}^{2}\left(\Omega\right)$
and $H_{q-1}(\nabla u)\in L_{loc}^{p}\left(\Omega\right)$, where
$p=q'=q/(q-1)$. Let $e_{s}$ be a coordinate direction in $\mathbb{R}^{n}$;
from now on, by a slight abuse of notation, we shall write $\tau_{h}g(x)$
to denote $g(x+he_{s})$, where $h\in\mathbb{R}\setminus\left\{ 0\right\} $.\\
$\hspace*{1em}$Now, let $\varphi\in W^{1,q}\left(\Omega\right)$
be compactly supported in $\Omega$ (i.e. a test function) and $h\in\mathbb{R}\setminus\left\{ 0\right\} $
be such that $\left|h\right|<\mathrm{dist}\left(\mathrm{supp}\left(\varphi\right),\mathbb{R}^{n}\setminus\Omega\right)$.
Since $u$ is a local weak solution of equation (\ref{eq:div}), we
have
\begin{equation}
\int_{\Omega}\langle\Delta_{h}H_{q-1}(\nabla u),\nabla\varphi\rangle\,dx=\int_{\Omega}\left(\Delta_{h}f\right)\varphi\,dx.\label{eq:start}
\end{equation}
\\
Let us fix two concentric balls $B_{R/2}$ and $B_{R}$ with $B_{R}\Subset\Omega$,
and consider a cut off function $\xi\in C_{c}^{\infty}\left(B_{R}\right)$
such that $\xi\equiv1$ on $\overline{B}_{R/2}$ and $\Vert\nabla\xi\Vert_{\infty}\leq C/R$.
Let $h\in\mathbb{R}\setminus\left\{ 0\right\} $ be such that $\left|h\right|\leq r_{0}<\frac{1}{2}\,\mathrm{dist}\left(B_{R},\mathbb{R}^{n}\setminus\Omega\right)$.
In what follows, we will denote by $c_{k}$ and $c$ some positive
constants which  do not depend on $h$, but they may vary on different
occasions.\\
$\hspace*{1em}$Now, let us consider, first, the case $f\in B_{p,\infty,loc}^{\alpha}\left(\Omega\right)$.
We then choose  $\varphi=\xi^{2}\Delta_{h}u$ as a test function into
(\ref{eq:start}). Using the fact that $u\in W_{loc}^{1,q}\left(\Omega\right)$,
$f\in B_{p,\infty,loc}^{\alpha}\left(\Omega\right)$, Hölder's inequality
and the properties of the difference quotients, and defining $B':=B_{R}+B(0,r_{0})$,
we get\begin{align*}
&\int_{\Omega}\langle\Delta_{h}H_{q-1}(\nabla u),\xi^{2}\Delta_{h}\nabla u+2\,\xi\nabla\xi\,\Delta_{h}u\rangle\,dx=\int_{\Omega}\langle\Delta_{h}H_{q-1}(\nabla u),\nabla\varphi\rangle\,dx\\
&=\int_{B_{R}}\left(\Delta_{h}f\right)\varphi\,dx\leq\int_{B_{R}}\left|\Delta_{h}f\right|\xi^{2}\left|\Delta_{h}u\right|\,dx\leq\left(\int_{B_{R}}\left|\Delta_{h}f\right|^{p}dx\right)^{\frac{1}{p}}\left(\int_{B_{R}}\xi^{2q}\left|\Delta_{h}u\right|^{q}dx\right)^{\frac{1}{q}}\\
&\leq\left(\int_{B_{R}}\frac{\left|\Delta_{h}f\right|^{p}}{\left|h\right|^{\alpha p}}\,dx\right)^{\frac{1}{p}}\left|h\right|^{\alpha}\,\Vert\Delta_{h}u\Vert_{L^{q}\left(B_{R}\right)}\leq c_{1}(n)\,\left(\int_{B_{R}}\frac{\left|\tau_{h}f-f\right|^{p}}{\left|h\right|^{\alpha p}}\,dx\right)^{\frac{1}{p}}\left|h\right|^{\alpha-1}\Vert\nabla u\Vert_{L^{q}\left(B'\right)}\\
&\leq c_{1}(n)\,\Vert f\Vert_{B_{p,\infty}^{\alpha}\left(B'\right)}\,\Vert\nabla u\Vert_{L^{q}\left(B'\right)}\left|h\right|^{\alpha-1}.
\end{align*} \\
From the previous estimate, we then obtain\begin{equation}\label{eq:stima0}
\begin{split}
&\int_{\Omega}\xi^{2}\langle\Delta_{h}H_{q-1}(\nabla u),\Delta_{h}\nabla u\rangle\,dx\\
&\leq c_{1}(n)\,\Vert f\Vert_{B_{p,\infty}^{\alpha}\left(B'\right)}\,\Vert\nabla u\Vert_{L^{q}\left(B'\right)}\left|h\right|^{\alpha-1}-2\int_{\Omega}\langle\Delta_{h}H_{q-1}(\nabla u),\nabla\xi\rangle\,\xi\,\Delta_{h}u\,dx.
\end{split}
\end{equation}Now, by Lemma \ref{lem:below} we have 
\begin{equation}
\beta(q,n)\int_{B_{R}}\xi^{2}\left|\Delta_{h}H_{\frac{q}{2}}(\nabla u(x))\right|^{2}dx\leq\int_{\Omega}\xi^{2}\langle\Delta_{h}H_{q-1}(\nabla u(x)),\Delta_{h}\nabla u(x)\rangle\,dx,\label{eq:stima01}
\end{equation}
where $\beta(q,n)$ is a positive constant. Combining estimates \eqref{eq:stima0}
and (\ref{eq:stima01}) and applying the Cauchy-Schwarz inequality
as well as the properties of $\xi$, we get \begin{align*}
&\beta(q,n)\int_{B_{R}}\xi^{2}\left|\Delta_{h}H_{\frac{q}{2}}(\nabla u)\right|^{2}dx\\
&\leq c_{1}(n)\,\Vert f\Vert_{B_{p,\infty}^{\alpha}\left(B'\right)}\,\Vert\nabla u\Vert_{L^{q}\left(B'\right)}\left|h\right|^{\alpha-1}-2\int_{\Omega}\langle\Delta_{h}H_{q-1}(\nabla u),\nabla\xi\rangle\,\xi\,\Delta_{h}u\,dx\\
&\leq c_{1}(n)\,\Vert f\Vert_{B_{p,\infty}^{\alpha}\left(B'\right)}\,\Vert\nabla u\Vert_{L^{q}\left(B'\right)}\left|h\right|^{\alpha-1}+\,c_{2}(R)\int_{\Omega}\left|\Delta_{h}H_{q-1}(\nabla u)\right|\xi\left|\Delta_{h}u\right|\,dx,
\end{align*}and dividing by $\beta(q,n)$, we obtain\begin{equation}\label{eq:STIMA1}
\begin{split}
&\int_{B_{R}}\xi^{2}\left|\Delta_{h}H_{\frac{q}{2}}(\nabla u)\right|^{2}dx\\
&\leq c_{3}(q,n)\,\Vert f\Vert_{B_{p,\infty}^{\alpha}\left(B'\right)}\,\Vert\nabla u\Vert_{L^{q}\left(B'\right)}\left|h\right|^{\alpha-1}+\,c_{4}(R,q,n)\int_{\Omega}\left|\Delta_{h}H_{q-1}(\nabla u)\right|\xi\left|\Delta_{h}u\right|\,dx.
\end{split}
\end{equation}Now we set 
\[
I_{1}:=\int_{\Omega}\left|\Delta_{h}H_{q-1}(\nabla u)\right|\xi\left|\Delta_{h}u\right|\,dx
\]
and apply Hölder's inequality with exponents $(q,p)$ and the properties
of the difference quotients to estimate $I_{1}$ as follows: 
\begin{equation}
I_{1}\leq c_{1}(n)\,\Vert\nabla u\Vert_{L^{q}\left(B'\right)}\left(\int_{B_{R}}\xi^{p}\left|\Delta_{h}H_{q-1}(\nabla u)\right|^{p}dx\right)^{\frac{1}{p}}.\label{eq:estI1}
\end{equation}
Similarly, we set 
\[
I_{2}:=\int_{B_{R}}\xi^{p}\left|\Delta_{h}H_{q-1}(\nabla u)\right|^{p}dx
\]
and use Lemma \ref{lem:jolly}, Hölder's inequality with exponents
$\left(\frac{2}{q},\frac{2}{2-q}\right)$ and the properties of the
difference quotients to control $I_{2}$ as follows:
\[
I_{2}\leq c_{5}(q,n)\,\int_{B_{R}}\xi^{p}\left|\Delta_{h}H_{1}(\nabla u)\right|^{p}\left[\left(\left|\tau_{h}\nabla u\right|-1\right)_{+}+\left(\left|\nabla u\right|-1\right)_{+}\right]^{p\left(q-2\right)}dx
\]
\[
=c_{5}(q,n)\,\int_{B_{R}}\xi^{p}\left|\Delta_{h}H_{1}(\nabla u)\right|^{q}\left|\Delta_{h}H_{1}(\nabla u)\right|^{p-q}\left[\left(\left|\tau_{h}\nabla u\right|-1\right)_{+}+\left(\left|\nabla u\right|-1\right)_{+}\right]^{p\left(q-2\right)}dx
\]
\[
=c_{5}(q,n)\left|h\right|^{q-p}\int_{B_{R}}\xi^{p}\left|\Delta_{h}H_{1}(\nabla u)\right|^{q}\left|H_{1}(\tau_{h}\nabla u)-H_{1}(\nabla u)\right|^{p-q}\left[\left(\left|\tau_{h}\nabla u\right|-1\right)_{+}+\left(\left|\nabla u\right|-1\right)_{+}\right]^{p\left(q-2\right)}dx
\]
\[
\leq c_{5}(q,n)\left|h\right|^{q-p}\int_{B_{R}}\xi^{p}\left|\Delta_{h}H_{1}(\nabla u)\right|^{q}\left[\left(\left|\tau_{h}\nabla u\right|-1\right)_{+}+\left(\left|\nabla u\right|-1\right)_{+}\right]^{pq-(p+q)}dx
\]
\[
=c_{5}(q,n)\left|h\right|^{q-p}\int_{B_{R}}\xi^{p}\left|\Delta_{h}H_{1}(\nabla u)\right|^{q}\left[\left(\left|\tau_{h}\nabla u\right|-1\right)_{+}+\left(\left|\nabla u\right|-1\right)_{+}\right]^{\frac{q}{2}\left(q-2\right)+\frac{q}{2}\left(2-q\right)}dx
\]
\[
\leq c_{5}(q,n)\left|h\right|^{q-p}\left(\int_{B_{R}}\xi^{\frac{2}{q-1}}\left|\Delta_{h}H_{1}(\nabla u)\right|^{2}\left[\left(\left|\tau_{h}\nabla u\right|-1\right)_{+}+\left(\left|\nabla u\right|-1\right)_{+}\right]^{q-2}dx\right)^{\frac{q}{2}}2^{\frac{q\left(2-q\right)}{2}}\,\Vert\nabla u\Vert_{L^{q}\left(B'\right)}^{\frac{q\left(2-q\right)}{2}}
\]
\begin{equation}
\leq c_{6}(q,n)\,\Vert\nabla u\Vert_{L^{q}\left(B'\right)}^{\frac{q\left(2-q\right)}{2}}\,\left|h\right|^{q-p}\left(\int_{B_{R}}\xi^{2}\left|\Delta_{h}H_{\frac{q}{2}}(\nabla u)\right|^{2}dx\right)^{\frac{q}{2}},\label{eq:est1}
\end{equation}
\\
where, in the last line, we have used Lemma \ref{lem:jolly} again,
as well as the properties of $\xi$. Collecting estimates \eqref{eq:STIMA1},
(\ref{eq:estI1}) and (\ref{eq:est1}) and applying Young's inequality
with $\theta>0$ and exponents $\left(\frac{2}{q-1},\frac{2}{3-q}\right)$,
we obtain\begin{equation}\label{eq:young1}
\begin{split}
&\int_{B_{R}}\xi^{2}\left|\Delta_{h}H_{\frac{q}{2}}(\nabla u)\right|^{2}dx\leq c_{3}(q,n)\,\Vert f\Vert_{B_{p,\infty}^{\alpha}\left(B'\right)}\,\Vert\nabla u\Vert_{L^{q}\left(B'\right)}\left|h\right|^{\alpha-1}+\,c_{7}(R,q,n)\,\Vert\nabla u\Vert_{L^{q}\left(B'\right)}\,I_{2}^{1/p}\\
&\leq c_{3}(q,n)\,\Vert f\Vert_{B_{p,\infty}^{\alpha}\left(B'\right)}\,\Vert\nabla u\Vert_{L^{q}\left(B'\right)}\left|h\right|^{\alpha-1}+\,c_{8}(R,q,n)\,\Vert\nabla u\Vert_{L^{q}\left(B'\right)}^{\frac{q(3-q)}{2}}\left|h\right|^{q-2}\left(\int_{B_{R}}\xi^{2}\left|\Delta_{h}H_{\frac{q}{2}}(\nabla u)\right|^{2}dx\right)^{\frac{q-1}{2}}\\
&\leq c_{3}(q,n)\,\Vert f\Vert_{B_{p,\infty}^{\alpha}\left(B'\right)}\,\Vert\nabla u\Vert_{L^{q}\left(B'\right)}\left|h\right|^{\alpha-1}\\
&\,\,\,\,+\,\frac{q-1}{2}\,\,\theta^{\frac{1}{q-1}}\int_{B_{R}}\xi^{2}\left|\Delta_{h}H_{\frac{q}{2}}(\nabla u)\right|^{2}dx\,+\,\frac{3-q}{2}\,\,\theta^{\frac{1}{q-3}}\,c_{9}(R,q,n)\,\Vert\nabla u\Vert_{L^{q}\left(B'\right)}^{q}\left|h\right|^{\frac{2(q-2)}{3-q}}.
\end{split}
\end{equation}Choosing $\theta=\left(\frac{1}{q-1}\right)^{q-1}$ and reabsorbing
the integral in the right-hand side of \eqref{eq:young1} by the left-hand
side, we get
\[
\int_{B_{R}}\xi^{2}\left|\Delta_{h}H_{\frac{q}{2}}(\nabla u)\right|^{2}dx\leq c\,\Vert f\Vert_{B_{p,\infty}^{\alpha}\left(B'\right)}\,\Vert\nabla u\Vert_{L^{q}\left(B'\right)}\left|h\right|^{\alpha-1}+\,c\,\Vert\nabla u\Vert_{L^{q}\left(B'\right)}^{q}\left|h\right|^{\frac{2(q-2)}{3-q}},
\]
from which we can infer
\begin{equation}
\int_{B_{R/2}}\left|\Delta_{h}H_{\frac{q}{2}}(\nabla u)\right|^{2}dx\leq c\,\Vert f\Vert_{B_{p,\infty}^{\alpha}\left(B'\right)}\,\Vert\nabla u\Vert_{L^{q}\left(B'\right)}\left|h\right|^{\alpha-1}+\,c\,\Vert\nabla u\Vert_{L^{q}\left(B'\right)}^{q}\left|h\right|^{\frac{2(q-2)}{3-q}},\label{eq:new1}
\end{equation}
with $c=c(R,q,n)>0$. Now, let $\lambda=\min\left\{ \frac{\alpha+1}{2},\frac{1}{3-q}\right\} $.
Dividing both sides of (\ref{eq:new1}) by $\left|h\right|^{2\lambda-2}$,
we then have
\begin{equation}
\underset{B_{R/2}}{\int}\left|\frac{\tau_{h}H_{q/2}(\nabla u)-H_{q/2}(\nabla u)}{\left|h\right|^{\lambda}}\right|^{2}dx\leq c\,\left(\Vert f\Vert_{B_{p,\infty}^{\alpha}\left(B'\right)}\,\Vert\nabla u\Vert_{L^{q}\left(B'\right)}\left|h\right|^{\alpha+1-2\lambda}+\,\Vert\nabla u\Vert_{L^{q}\left(B'\right)}^{q}\left|h\right|^{\frac{2}{3-q}-2\lambda}\right).\label{eq:besvar3}
\end{equation}
Since the above estimate holds for every $h\in\mathbb{R}\setminus\left\{ 0\right\} $
such that $\left|h\right|\leq r_{0}$, we can take the supremum over
$\left|h\right|<\delta$ for some $\delta<r_{0}$ and obtain
\[
\sup_{\left|h\right|<\delta}\underset{B_{R/2}}{\int}\left|\frac{\tau_{h}H_{q/2}(\nabla u)-H_{q/2}(\nabla u)}{\left|h\right|^{\lambda}}\right|^{2}dx\leq c\,\left(\Vert f\Vert_{B_{p,\infty}^{\alpha}\left(B'\right)}\,\Vert\nabla u\Vert_{L^{q}\left(B'\right)}\,r_{0}^{\alpha+1-2\lambda}+\,\Vert\nabla u\Vert_{L^{q}\left(B'\right)}^{q}\,r_{0}^{\frac{2}{3-q}-2\lambda}\right).
\]
In particular, this means that $H_{\frac{q}{2}}(\nabla u)\in B_{2,\infty}^{\min\left\{ \frac{\alpha+1}{2},\frac{1}{3-q}\right\} }$
locally in $\Omega$. \\
$\hspace*{1em}$Finally, let us consider the case $f\in W_{loc}^{1,p}\left(\Omega\right)$.
Arguing as above, but this time using Lemma \ref{lem:diffquo} to
estimate the $L^{p}$-norm of the difference quotient of $f$, we
obtain 
\begin{equation}
\underset{B_{R/2}}{\int}\left|\frac{\tau_{h}H_{q/2}(\nabla u)-H_{q/2}(\nabla u)}{\left|h\right|^{1/(3-q)}}\right|^{2}dx\leq c\,\left(\Vert\nabla f\Vert_{L^{p}\left(B'\right)}\,\Vert\nabla u\Vert_{L^{q}\left(B'\right)}\left|h\right|^{\frac{2(2-q)}{3-q}}+\,\Vert\nabla u\Vert_{L^{q}\left(B'\right)}^{q}\right),\label{eq:besvar2}
\end{equation}
which holds for every $h\in\mathbb{R}\setminus\left\{ 0\right\} $
such that $\left|h\right|\leq r_{0}$. Arguing exactly as before,
we may then conclude that $H_{\frac{q}{2}}(\nabla u)\in B_{2,\infty,loc}^{\frac{1}{3-q}}\left(\Omega\right)$.\end{proof}

\noindent $\hspace*{1em}$We are now in position to give the

\noindent \begin{proof}[\bfseries{Proof of Theorem~\ref{thm:main2}}]
Let us begin with the case $f\in W_{loc}^{1,p}\left(\Omega\right)$.
By duality, we know that $\sigma_{0}$ is related to any solution
$u$ of problem (\ref{eq:P2}) through the optimality condition
\[
\sigma_{0}(x)=\nabla\mathcal{H}^{*}(\nabla u(x))=H_{q-1}(\nabla u(x)),\,\,\,\,\,\mathrm{for}\,\,\mathscr{L}^{n}\textrm{-}\mathrm{a.e.}\,\,x\in\Omega.
\]
Now we observe that for every $q\in(1,2)$ we have
\[
0<q-1<\frac{q}{2}.
\]
If $B_{R/2}$ is the same ball fixed in the proof of Theorem \ref{thm:mainth},
then applying Lemma \ref{lem:jolly} we get\begin{equation}\label{eq:besvar}
\begin{split}
&\underset{B_{R/2}}{\int}\left|\frac{\tau_{h}\sigma_{0}-\sigma_{0}}{\left|h\right|^{\frac{2}{p(3-q)}}}\right|^{p}dx\leq c(q,n)\left|h\right|^{\frac{2}{q-3}}\underset{B_{R/2}}{\int}\left|\tau_{h}H_{\frac{q}{2}}(\nabla u)-H_{\frac{q}{2}}(\nabla u)\right|^{p}\left[\left(\left|\tau_{h}\nabla u\right|-1\right)_{+}^{\frac{q}{2}}+\left(\left|\nabla u\right|-1\right)_{+}^{\frac{q}{2}}\right]^{\frac{q-2}{q-1}}dx\\
&=c(q,n)\underset{B_{R/2}}{\int}\left|\frac{\tau_{h}H_{q/2}(\nabla u)-H_{q/2}(\nabla u)}{\left|h\right|^{1/(3-q)}}\right|^{2}\left|\tau_{h}H_{\frac{q}{2}}(\nabla u)-H_{\frac{q}{2}}(\nabla u)\right|^{p-2}\left[\left(\left|\tau_{h}\nabla u\right|-1\right)_{+}^{\frac{q}{2}}+\left(\left|\nabla u\right|-1\right)_{+}^{\frac{q}{2}}\right]^{\frac{q-2}{q-1}}dx\\
&\leq c(q,n)\underset{B_{R/2}}{\int}\left|\frac{\tau_{h}H_{q/2}(\nabla u)-H_{q/2}(\nabla u)}{\left|h\right|^{1/(3-q)}}\right|^{2}dx.
\end{split}
\end{equation}Combining estimates \eqref{eq:besvar} and (\ref{eq:besvar2}) and
arguing as in the last part of the proof of Theorem \ref{thm:mainth},
we obtain the desired conclusion.\\
$\hspace*{1em}$Finally, when $f\in B_{p,\infty,loc}^{\alpha}\left(\Omega\right)$,
a similar argument applies, but taking into account estimate (\ref{eq:besvar3})
instead of (\ref{eq:besvar2}).\end{proof}

\noindent $\hspace*{1em}$Now we observe that the Besov regularity
of $\sigma_{0}$ established in Theorem  \ref{thm:main2} allows us
to get a gain of integrability for $\sigma_{0}$. More precisely,
we have:

\noindent \begin{proof}[\bfseries{Proof of Corollary~\ref{cor:corol1}}]
We first note that for every $q\in(1,2)$ and $\alpha\in(0,1)$ we
have 
\[
\frac{2}{3-q}<2\leq n\,\,\,\,\,\,\,\,\mathrm{and}\,\,\,\,\,\,\,\,\alpha+1<n.
\]
$\hspace*{1em}$If $f\in B_{p,\infty,loc}^{\alpha}\left(\Omega\right)$
and $\frac{\alpha+1}{p}\geq\frac{2}{p(3-q)}$ (that is, $\frac{q-1}{3-q}\leq\alpha<1$),
then  from Theorems \ref{thm:main2} and \ref{thm:emb2} it follows
that $\sigma_{0}\in L_{loc}^{s}\left(\Omega\right)$ for all $s\in\left[1,\frac{np(3-q)}{n(3-q)-2}\right)$.\\
$\hspace*{1em}$When $f\in W_{loc}^{1,p}\left(\Omega\right)$, thanks
to Theorem \ref{thm:main2} we have that $\sigma_{0}\in B_{p,\infty,loc}^{\frac{2}{p(3-q)}}\left(\Omega\right)$.
We therefore reach the same conclusion as in the previous case. \\
$\hspace*{1em}$Finally, when $f\in B_{p,\infty,loc}^{\alpha}\left(\Omega\right)$
and $\frac{\alpha+1}{p}<\frac{2}{p(3-q)}$ (that is, $0<\alpha<\frac{q-1}{3-q}$),
Theorems \ref{thm:main2} and \ref{thm:emb2} imply that $\sigma_{0}\in L_{loc}^{s}\left(\Omega\right)$
for all $s\in\left[1,\frac{np}{n-\alpha-1}\right)$.\end{proof} 

\noindent $\hspace*{1em}$As an easy consequence of the (local) higher
integrability of $H_{q-1}(\nabla u)$ established in the proof of
Corollary  \ref{cor:corol1}, we get a gain of integrability for $\nabla u$.
More precisely, we have:

\noindent \begin{proof}[\bfseries{Proof of Corollary~\ref{cor:corol3}}]
In the case $f\in W_{loc}^{1,p}\left(\Omega\right)$, we get $H_{q-1}(\nabla u)\in L_{loc}^{s}\left(\Omega\right)$
for every $s\in\left[1,\frac{np(3-q)}{n(3-q)-2}\right)$ and then
\[
\int_{K}\left(\left|\nabla u(x)\right|-1\right)_{+}^{(q-1)s}dx=\int_{K}\left|H_{q-1}(\nabla u(x))\right|^{s}dx<+\infty
\]
for all $s\in\left[1,\frac{np(3-q)}{n(3-q)-2}\right)$ and all compact
subsets $K$ of $\Omega$, which ensures that $\nabla u\in L_{loc}^{r}\left(\Omega\right)$
for all  $r\in\left[1,\frac{nq(3-q)}{n(3-q)-2}\right)$. In the other
two cases, a similar argument applies, thus proving the assertion.\end{proof}
\end{singlespace}
\begin{singlespace}

\section{The non-singular case $q\protect\geq2$\label{sec:q>2}}
\end{singlespace}

\begin{singlespace}
\noindent $\hspace*{1em}$Here we proceed with the analysis of the
non-singular case, by proving that, unlike what may occur in the singular
one, the Besov regularity of the datum $f$ translates into a Besov
regularity for $H_{\frac{q}{2}}(\nabla u)$ with no loss in the order
of differentiation. More precisely, we have the following\\

\end{singlespace}
\begin{thm}
\begin{singlespace}
\noindent \label{thm:q>2} Let $n\geq2$, $q\geq2$, $\alpha\in\left(0,1\right)$
and $f\in B_{p,\infty,loc}^{\alpha}\left(\Omega\right)$. Moreover,
let $u\in W_{loc}^{1,q}\left(\Omega\right)$ be a local weak solution
of the equation
\begin{equation}
\mathrm{-div}\,(\nabla\mathcal{H}^{*}(\nabla u))=f.\label{eq:div1}
\end{equation}
Then $H_{\frac{q}{2}}(\nabla u)\in B_{2,\infty,loc}^{\frac{\alpha+1}{2}}\left(\Omega\right)$.
Furthermore, for any ball $B_{R}\Subset\Omega$, the following estimate
\[
\underset{B_{R/2}}{\int}\left|\tau_{s,h}H_{\frac{q}{2}}(\nabla u)-H_{\frac{q}{2}}(\nabla u)\right|^{2}dx\leq c\,\left(\Vert f\Vert_{B_{p,\infty}^{\alpha}\left(B'\right)}\,\Vert\nabla u\Vert_{L^{q}\left(B'\right)}\left|h\right|^{\alpha+1}+\,\Vert\nabla u\Vert_{L^{q}\left(B'\right)}^{q}\left|h\right|^{2}\right)
\]
holds true for every $s\in\left\{ 1,\ldots,n\right\} $, for every
$h\in\mathbb{R}$ such that $\left|h\right|\leq r_{0}<\frac{1}{2}\,\mathrm{dist}\left(B_{R},\partial\Omega\right)$,
for $B'=B_{R}+B(0,r_{0})$ and a positive constant $c=c(R,q,n)$.
\end{singlespace}
\end{thm}

\begin{singlespace}
\noindent \begin{proof}[\bfseries{Proof}] We first observe that $H_{\frac{q}{2}}(\nabla u)\in L_{loc}^{2}\left(\Omega\right)$
and $H_{q-1}(\nabla u)\in L_{loc}^{p}\left(\Omega\right)$. Let $e_{s}$
be a coordinate direction in $\mathbb{R}^{n}$; from now on, by a
slight abuse of notation, we shall write $\tau_{h}g(x)$ to denote
$g(x+he_{s})$, where $h\in\mathbb{R}\setminus\left\{ 0\right\} $.\\
$\hspace*{1em}$Now, let $\varphi\in W^{1,q}\left(\Omega\right)$
be compactly supported in $\Omega$ (i.e. a test function) and $h\in\mathbb{R}\setminus\left\{ 0\right\} $
be such that $\left|h\right|<\mathrm{dist}\left(\mathrm{supp}\left(\varphi\right),\mathbb{R}^{n}\setminus\Omega\right)$.
Since $u$ is a local weak solution of equation (\ref{eq:div1}),
we have
\begin{equation}
\int_{\Omega}\langle\Delta_{h}H_{q-1}(\nabla u),\nabla\varphi\rangle\,dx=\int_{\Omega}\left(\Delta_{h}f\right)\varphi\,dx.\label{eq:start-1}
\end{equation}
\\
Let us fix two concentric balls $B_{R/2}$ and $B_{R}$ with $B_{R}\Subset\Omega$,
and consider a cut off function $\xi\in C_{c}^{\infty}\left(B_{R}\right)$
such that $\xi\equiv1$ on $\overline{B}_{R/2}$ and $\Vert\nabla\xi\Vert_{\infty}\leq C/R$.
Let $h\in\mathbb{R}\setminus\left\{ 0\right\} $ be such that $\left|h\right|\leq r_{0}<\frac{1}{2}\,\mathrm{dist}\left(B_{R},\mathbb{R}^{n}\setminus\Omega\right)$.
In what follows, we will denote by $c_{k}$ and $c$ some positive
constants which  do not depend on $h$, but they may vary on different
occasions. We then choose $\varphi=\xi^{2}\Delta_{h}u$ as a test
function into (\ref{eq:start-1}). Using the fact that $u\in W_{loc}^{1,q}\left(\Omega\right)$,
$f\in B_{p,\infty,loc}^{\alpha}\left(\Omega\right)$, Hölder's inequality
and the properties of the difference quotients, and defining $B':=B_{R}+B(0,r_{0})$,
we get\begin{align*}
&\int_{\Omega}\langle\Delta_{h}H_{q-1}(\nabla u),\xi^{2}\Delta_{h}\nabla u+2\,\xi\nabla\xi\,\Delta_{h}u\rangle\,dx=\int_{\Omega}\langle\Delta_{h}H_{q-1}(\nabla u),\nabla\varphi\rangle\,dx\\
&=\int_{B_{R}}\left(\Delta_{h}f\right)\varphi\,dx\leq\int_{B_{R}}\left|\Delta_{h}f\right|\xi^{2}\left|\Delta_{h}u\right|\,dx\leq\left(\int_{B_{R}}\left|\Delta_{h}f\right|^{p}dx\right)^{\frac{1}{p}}\left(\int_{B_{R}}\xi^{2q}\left|\Delta_{h}u\right|^{q}dx\right)^{\frac{1}{q}}\\
&\leq\left(\int_{B_{R}}\frac{\left|\Delta_{h}f\right|^{p}}{\left|h\right|^{\alpha p}}\,dx\right)^{\frac{1}{p}}\left|h\right|^{\alpha}\,\Vert\Delta_{h}u\Vert_{L^{q}\left(B_{R}\right)}\leq c_{1}(n)\,\left(\int_{B_{R}}\frac{\left|\tau_{h}f-f\right|^{p}}{\left|h\right|^{\alpha p}}\,dx\right)^{\frac{1}{p}}\left|h\right|^{\alpha-1}\Vert\nabla u\Vert_{L^{q}\left(B'\right)}\\
&\leq c_{1}(n)\,\Vert f\Vert_{B_{p,\infty}^{\alpha}\left(B'\right)}\,\Vert\nabla u\Vert_{L^{q}\left(B'\right)}\left|h\right|^{\alpha-1},
\end{align*}from which we obtain\begin{equation}\label{eq:stima0-1}
\begin{split}
&\int_{\Omega}\xi^{2}\langle\Delta_{h}H_{q-1}(\nabla u),\Delta_{h}\nabla u\rangle\,dx\\
&\leq c_{1}(n)\,\Vert f\Vert_{B_{p,\infty}^{\alpha}\left(B'\right)}\,\Vert\nabla u\Vert_{L^{q}\left(B'\right)}\left|h\right|^{\alpha-1}-2\int_{\Omega}\langle\Delta_{h}H_{q-1}(\nabla u),\nabla\xi\rangle\,\xi\,\Delta_{h}u\,dx.
\end{split}
\end{equation}Now, by Lemma \ref{lem:Brasco} we have 
\[
\frac{4}{q^{2}}\int_{B_{R}}\xi^{2}\left|\Delta_{h}H_{\frac{q}{2}}(\nabla u(x))\right|^{2}dx\leq\int_{\Omega}\xi^{2}\langle\Delta_{h}H_{q-1}(\nabla u(x)),\Delta_{h}\nabla u(x)\rangle\,dx.
\]
\\
Using the above estimate together with \eqref{eq:stima0-1} and the
Cauchy-Schwarz inequality, we then get\begin{align*}
&\frac{4}{q^{2}}\int_{B_{R}}\xi^{2}\left|\Delta_{h}H_{\frac{q}{2}}(\nabla u)\right|^{2}dx\\
&\leq c_{1}(n)\,\Vert f\Vert_{B_{p,\infty}^{\alpha}\left(B'\right)}\,\Vert\nabla u\Vert_{L^{q}\left(B'\right)}\left|h\right|^{\alpha-1}-2\int_{\Omega}\langle\Delta_{h}H_{q-1}(\nabla u),\nabla\xi\rangle\,\xi\,\Delta_{h}u\,dx\\
&\leq c_{1}(n)\,\Vert f\Vert_{B_{p,\infty}^{\alpha}\left(B'\right)}\,\Vert\nabla u\Vert_{L^{q}\left(B'\right)}\left|h\right|^{\alpha-1}+\,2\,\Vert\nabla\xi\Vert_{L^{\infty}}\int_{\Omega}\left|\Delta_{h}H_{q-1}(\nabla u)\right|\xi\left|\Delta_{h}u\right|\,dx\\
&\leq c_{1}(n)\,\Vert f\Vert_{B_{p,\infty}^{\alpha}\left(B'\right)}\,\Vert\nabla u\Vert_{L^{q}\left(B'\right)}\left|h\right|^{\alpha-1}+\,c_{2}(R)\int_{\Omega}\left|\Delta_{h}H_{q-1}(\nabla u)\right|\xi\left|\Delta_{h}u\right|\,dx,
\end{align*}and dividing by $4/q^{2}$, we obtain\begin{equation}\label{eq:STIMA1-1}
\begin{split}
&\int_{B_{R}}\xi^{2}\left|\Delta_{h}H_{\frac{q}{2}}(\nabla u)\right|^{2}dx\\
&\leq c_{3}(q,n)\,\Vert f\Vert_{B_{p,\infty}^{\alpha}\left(B'\right)}\,\Vert\nabla u\Vert_{L^{q}\left(B'\right)}\left|h\right|^{\alpha-1}+\,c_{4}(R,q)\int_{\Omega}\left|\Delta_{h}H_{q-1}(\nabla u)\right|\xi\left|\Delta_{h}u\right|\,dx.
\end{split}
\end{equation}Now, our aim is to estimate 
\[
I_{1}:=\int_{\Omega}\left|\Delta_{h}H_{q-1}(\nabla u)\right|\xi\left|\Delta_{h}u\right|\,dx.
\]
By virtue of Lemma \ref{lem:Brasco} we have
\begin{equation}
\left|\Delta_{h}H_{q-1}(\nabla u(x))\right|\leq(q-1)\left(\left|\tau_{h}H_{\frac{q}{2}}(\nabla u(x))\right|^{\frac{q-2}{q}}+\left|H_{\frac{q}{2}}(\nabla u(x))\right|^{\frac{q-2}{q}}\right)\left|\Delta_{h}H_{\frac{q}{2}}(\nabla u(x))\right|.\label{eq:jo1-1}
\end{equation}
We now use (\ref{eq:jo1-1}), Hölder's inequality with exponents $\left(q,2,\frac{2q}{q-2}\right)$
and the properties of the difference quotients to control $I_{1}$
as follows:\\
 \begin{equation}\label{eq:STIMA11}
\begin{split}
I_{1}&\leq(q-1)\int_{B_{R}}\xi\left|\Delta_{h}u\right|\left(\left|\tau_{h}H_{\frac{q}{2}}(\nabla u)\right|^{\frac{q-2}{q}}+\left|H_{\frac{q}{2}}(\nabla u)\right|^{\frac{q-2}{q}}\right)\left|\Delta_{h}H_{\frac{q}{2}}(\nabla u)\right|dx\\
             &\leq(q-1)\,\Vert\Delta_{h}u\Vert_{L^{q}\left(B_{R}\right)}\left(\int_{B_{R}}\xi^{2}\left|\Delta_{h}H_{\frac{q}{2}}(\nabla u)\right|^{2}dx\right)^{\frac{1}{2}}\\
&\,\,\,\,\cdot\left(\int_{B_{R}}\left(\left|\tau_{h}H_{\frac{q}{2}}(\nabla u)\right|^{\frac{q-2}{q}}+\left|H_{\frac{q}{2}}(\nabla u)\right|^{\frac{q-2}{q}}\right)^{\frac{2q}{q-2}}dx\right)^{\frac{q-2}{2q}}\\
&\leq c_{5}(q,n)\,\Vert\nabla u\Vert_{L^{q}\left(B'\right)}\left(\int_{B_{R}}\xi^{2}\left|\Delta_{h}H_{\frac{q}{2}}(\nabla u)\right|^{2}dx\right)^{\frac{1}{2}}2\,\big\Vert H_{\frac{q}{2}}(\nabla u)\big\Vert_{L^{2}\left(B'\right)}^{\frac{q-2}{q}}\\
&\leq c_{6}(q,n)\,\Vert\nabla u\Vert_{L^{q}\left(B'\right)}^{q/2}\left(\int_{B_{R}}\xi^{2}\left|\Delta_{h}H_{\frac{q}{2}}(\nabla u)\right|^{2}dx\right)^{\frac{1}{2}}.
\end{split}
\end{equation}\\
Joining \eqref{eq:STIMA1-1} and \eqref{eq:STIMA11} and applying
Young's inequality with exponents $\left(2,2\right)$, we obtain\begin{equation}\label{eq:young2}
\begin{split}
&\int_{B_{R}}\xi^{2}\left|\Delta_{h}H_{\frac{q}{2}}(\nabla u)\right|^{2}dx\leq c_{3}(q,n)\,\Vert f\Vert_{B_{p,\infty}^{\alpha}\left(B'\right)}\,\Vert\nabla u\Vert_{L^{q}\left(B'\right)}\left|h\right|^{\alpha-1}+\,c_{4}(R,q)\,I_{1}\\
&\leq c_{3}(q,n)\,\Vert f\Vert_{B_{p,\infty}^{\alpha}\left(B'\right)}\,\Vert\nabla u\Vert_{L^{q}\left(B'\right)}\left|h\right|^{\alpha-1}+\,c_{7}(R,q,n)\,\Vert\nabla u\Vert_{L^{q}\left(B'\right)}^{q/2}\left(\int_{B_{R}}\xi^{2}\left|\Delta_{h}H_{\frac{q}{2}}(\nabla u)\right|^{2}dx\right)^{\frac{1}{2}}\\
&\leq c_{3}(q,n)\,\Vert f\Vert_{B_{p,\infty}^{\alpha}\left(B'\right)}\,\Vert\nabla u\Vert_{L^{q}\left(B'\right)}\left|h\right|^{\alpha-1}+\,\frac{1}{2}\,c_{8}(R,q,n)\,\Vert\nabla u\Vert_{L^{q}\left(B'\right)}^{q}+\,\frac{1}{2}\int_{B_{R}}\xi^{2}\left|\Delta_{h}H_{\frac{q}{2}}(\nabla u)\right|^{2}dx.
\end{split}
\end{equation}Reabsorbing the integral in the right-hand side of \eqref{eq:young2}
by the left-hand side, we get
\[
\int_{B_{R}}\xi^{2}\left|\Delta_{h}H_{\frac{q}{2}}(\nabla u)\right|^{2}dx\leq c\,\left(\Vert f\Vert_{B_{p,\infty}^{\alpha}\left(B'\right)}\,\Vert\nabla u\Vert_{L^{q}\left(B'\right)}\left|h\right|^{\alpha-1}+\,\Vert\nabla u\Vert_{L^{q}\left(B'\right)}^{q}\right),
\]
from which we can infer
\begin{equation}
\int_{B_{R/2}}\left|\Delta_{h}H_{\frac{q}{2}}(\nabla u)\right|^{2}dx\leq c\,\left(\Vert f\Vert_{B_{p,\infty}^{\alpha}\left(B'\right)}\,\Vert\nabla u\Vert_{L^{q}\left(B'\right)}\left|h\right|^{\alpha-1}+\,\Vert\nabla u\Vert_{L^{q}\left(B'\right)}^{q}\right),\label{eq:new1-1}
\end{equation}
with $c=c(R,q,n)>0$. Multiplying both sides of (\ref{eq:new1-1})
by $\left|h\right|^{1-\alpha}$, we then have
\[
\underset{B_{R/2}}{\int}\left|\frac{\tau_{h}H_{q/2}(\nabla u)-H_{q/2}(\nabla u)}{\left|h\right|^{(\alpha+1)/2}}\right|^{2}dx\leq c\,\left(\Vert f\Vert_{B_{p,\infty}^{\alpha}\left(B'\right)}\,\Vert\nabla u\Vert_{L^{q}\left(B'\right)}\,+\,\Vert\nabla u\Vert_{L^{q}\left(B'\right)}^{q}\left|h\right|^{1-\alpha}\right).
\]
Since the above estimate holds for every $h\in\mathbb{R}\setminus\left\{ 0\right\} $
such that $\left|h\right|\leq r_{0}$, we can take the supremum over
$\left|h\right|<\delta$ for some $\delta<r_{0}$ and obtain
\[
\sup_{\left|h\right|<\delta}\underset{B_{R/2}}{\int}\left|\frac{\tau_{h}H_{q/2}(\nabla u)-H_{q/2}(\nabla u)}{\left|h\right|^{(\alpha+1)/2}}\right|^{2}dx\leq c\,\left(\Vert f\Vert_{B_{p,\infty}^{\alpha}\left(B'\right)}\,\Vert\nabla u\Vert_{L^{q}\left(B'\right)}\,+\,\Vert\nabla u\Vert_{L^{q}\left(B'\right)}^{q}\,r_{0}^{1-\alpha}\right),
\]
which gives the desired conclusion.\end{proof}

\noindent $\hspace*{1em}$The Besov regularity of $H_{\frac{q}{2}}(\nabla u)$
established in the previous theorem allows us to get a gain of integrability
for $\nabla u$. More precisely, we have the following
\end{singlespace}
\begin{cor}
\begin{singlespace}
\noindent \label{cor:Cor1} Under the assumptions of Theorem \ref{thm:q>2},
we get 
\[
H_{\frac{q}{2}}(\nabla u)\in L_{loc}^{r}\left(\Omega\right)\,\,\,for\,\,all\,\,r\in\left[1,\frac{2n}{n-\alpha-1}\right)
\]
and
\[
\nabla u\in L_{loc}^{s}\left(\Omega\right)\,\,\,for\,\,all\,\,s\in\left[1,\frac{nq}{n-\alpha-1}\right).
\]
\end{singlespace}
\end{cor}

\begin{singlespace}
\noindent \begin{proof}[\bfseries{Proof}] By observing that $\alpha+1<n$,
from Theorems \ref{thm:q>2} and \ref{thm:emb2} we immediately obtain
the first conclusion. Then 
\[
\int_{K}\left(\left|\nabla u(x)\right|-1\right)_{+}^{\frac{qr}{2}}dx=\int_{K}\left|H_{\frac{q}{2}}(\nabla u(x))\right|^{r}dx<+\infty
\]
for all $r\in\left[1,\frac{2n}{n-\alpha-1}\right)$ and all compact
subsets $K$ of $\Omega$, which ensures that $\nabla u\in L_{loc}^{s}\left(\Omega\right)$
for all $s\in\left[1,\frac{nq}{n-\alpha-1}\right)$.\end{proof}

\noindent \begin{brem} Let us observe that when $n=2$ (the case
which is relevant for applications to network congestion), the previous
result implies that $\nabla u\in L_{loc}^{s}\left(\Omega\right)$
for all $s<2q/(1-\alpha)$. This means that for every $r>2q$ there
exists a $\delta=\delta(r,q)\in\left(0,1\right)$ such that, if $\alpha\in\left(\delta,1\right)$
and $f\in B_{p,\infty,loc}^{\alpha}\left(\Omega\right)$, then $\nabla u\in L_{loc}^{r}\left(\Omega\right)$.
In fact, it is enough to take $\delta=1-2q/r.$ In other words, given
any $r>2q$, we can choose $\alpha$ close enough to $1$ in order
to have $\nabla u\in L_{loc}^{r}\left(\Omega\right)$ whenever $f\in B_{p,\infty,loc}^{\alpha}\left(\Omega\right)$.\end{brem}\bigskip{}

\noindent $\hspace*{1em}$Again, if we come back to the variational
problem (\ref{eq:P1}), we obtain the following higher integrability
result for its minimizer:
\end{singlespace}
\begin{cor}
\begin{singlespace}
\noindent Let $n\geq2$, $q\geq2$, $\alpha\in\left(0,1\right)$ and
$f\in B_{p,\infty,loc}^{\alpha}\left(\Omega\right)$. Moreover, let
$\sigma_{0}\in L^{p}\left(\Omega,\mathbb{R}^{n}\right)$ be the solution
of $\left(\ref{eq:P1}\right)$. Then $\sigma_{0}\in L_{loc}^{r}\left(\Omega\right)$
for all $r\in\left[1,\frac{np}{n-\alpha-1}\right)$. 
\end{singlespace}
\end{cor}

\begin{singlespace}
\noindent \begin{proof}[\bfseries{Proof}] By duality, we know that
$\sigma_{0}$ is related to any solution $u$ of problem (\ref{eq:P2})
through the optimality condition
\[
\sigma_{0}(x)=\nabla\mathcal{H}^{*}(\nabla u(x))=\left(\left|\nabla u(x)\right|-1\right)_{+}^{q-1}\frac{\nabla u(x)}{\left|\nabla u(x)\right|},\,\,\,\,\,\mathrm{for}\,\,\mathscr{L}^{n}\textrm{-}\mathrm{a.e.}\,\,x\in\Omega.
\]
Since $u\in W^{1,q}\left(\Omega\right)$ is a weak solution of the
Euler-Lagrange equation (\ref{eq:Eul-Lag}) and $\left|\sigma_{0}\right|=\left|H_{\frac{q}{2}}(\nabla u)\right|^{\frac{2}{p}}$
$\mathscr{L}^{n}$-a.e. in $\Omega$, the assertion  immediately follows
from Corollary \ref{cor:Cor1}.\end{proof}

\noindent \begin{brem} Let us observe that when $n=2$,  the previous
result implies that $\sigma_{0}\in L_{loc}^{r}\left(\Omega\right)$
for all $r<2p/(1-\alpha)$. This means that for every $s>2p$ there
exists a $\delta=\delta(s,p)\in\left(0,1\right)$ such that, if $\alpha\in\left(\delta,1\right)$
and $f\in B_{p,\infty,loc}^{\alpha}\left(\Omega\right)$, then $\sigma_{0}\in L_{loc}^{s}\left(\Omega\right)$.
In fact, it is enough to take $\delta=1-2p/s.$ In other words, given
any $s>2p$, we can choose $\alpha$ close enough to $1$ in order
to have $\sigma_{0}\in L_{loc}^{s}\left(\Omega\right)$ whenever $f\in B_{p,\infty,loc}^{\alpha}\left(\Omega\right)$.\end{brem}\bigskip{}

\noindent $\hspace*{1em}$\textbf{Acknowledgements.} The author gratefully
acknowledges fruitful discussions with Antonia Passarelli di Napoli
during the preparation of this paper.\\
$\hspace*{1em}$Moreover, he would like to thank the reviewer for
his/her suggestions to improve this work. 
\end{singlespace}

\begin{singlespace}

\lyxaddress{\noindent \textbf{$\quad$}\\
$\hspace*{1em}$\textbf{Pasquale Ambrosio}\\
Dipartimento di Matematica e Applicazioni ``R. Caccioppoli''\\
Università degli Studi di Napoli ``Federico II''\\
Via Cintia, 80126 Napoli, Italy.\\
\textit{E-mail address}: pasquale.ambrosio2@studenti.unina.it}
\end{singlespace}

\end{document}